\documentclass[draft]{amsart}

\theoremstyle{definition}

\theoremstyle{remark}

\numberwithin{equation}{section}

\usepackage{amsmath}
\usepackage{amsfonts}
\usepackage{amssymb}
\usepackage{graphicx}
\usepackage{tikz}
\usetikzlibrary{arrows.meta}
\usetikzlibrary{plotmarks}

\newcommand{\EE}{\mathbb{E}}

\begin{document}


\title{A model for runway operation decisions with stable queues}


\author{Carlos D.F.J. Bernardes}
\address{Escola Superior de Tecnologia de Setúbal, Instituto Politécnico de Setúbal}
\email{carlos.fidalgo2012@gmail.com}


\author{César Rodrigo}
\address{Escola Superior de Tecnologia de Setúbal, Instituto Politécnico de Setúbal \\ Centro de Matemática e Aplicações Fundamentais -- Centro de Investigação Operacional (CMAF-CIO)}
\email{cesar.fernandez@estsetubal.ips.pt}

\maketitle 

\begin{abstract}
The landing and takeoff operations for an airport at any given day are described in terms of the capacity envelopes associated to runway system configurations, of the scheduled flights along the day and of predefined delay tolerances for both types of operations. Assuming the inter-arrival times and service times are random variables with known quadratic ratio of momenta, it is possible to identify a parameter, the stable transit time associated to each service and slot, that measures the performance of the airport along the day. Even though constraints on service rates and definition of transit times are nonlinear, the description of runway system capacity using observed operational throughput control points allows to perform an optimization of the service given by the runway system in the general linear programming framework, minimizing costs associated to delays.
\end{abstract}

\section{Introduction}

Since the beginning of commercial aviation, nations and territories saw their social and industrial development linked to the degree of connectivity among them  that could be achieved by this means of transport, which allowed the expedition of passengers and commodities between cities and production centers in a more immediate way than using any of the means available before. The economical growth of cities has since largely depended on the level of service made possible by air traffic.

The total throughput that this industry could achieve was initially limited by technical capacities of existing airplanes, which progressively grew in both distance and load capabilities. Since the beginning of the 21st century this industry finds itself in a situation where the main factors that determine throughput limitations are not so much the technical specifications of airplanes but rather an adequate management of airport facilities, which strongly restrict the number of flights at the local airspace of a big city.

An appropriate management of airports seeks for a balance between, on the one hand, several parameters measuring productivity (in terms of transported cargo, monetary profit, timing, or environmental issues) and on the other hand restraints given by external regulations which aim to enforce interests, rights and security for all parties involved in the transportation process.

A classification of air traffic management problems can be found in \cite{Research_review_of_air_traffic_management}, and its treatment involves a large variety of mathematical techniques presented in a rich literature of operations research studies. An increasing fraction of these studies is devoted to management of operations within an airport's premises, as described in \cite{Planning_and_design_of_airports_fifth_edition}. An airport's runway system represents a bottleneck that limits the execution of services that demand runway occupancy and is critical for the productivity of the airport. It is an essential component that limits the airport throughput capabilities in landing (arrival) or takeoff (departure). A seminal work in this area which has influenced all subsequent studies was presented in 1960 by Blumstein \cite{Blumstein_AnaliticalInvestigationAirportCapacity}, containing an analytical model for the operational capacity of a runway as a function of relevant parameters regarding the airport and air traffic (security rules, meteorological conditions, traffic intensity, navigational instrumentation). The performance of the airport depends thus on the airport runway capacity, which can be modified only by large infrastructure investment and is assumed to be fixed. It also depends on the flight scheduling, a forecast of the services to be performed each day, which in many cases depends on decisions that are external to the local air traffic management, as it involves several airports and flight connections for different operators. Finally, the performance can be adjusted locally at the airport by the aplication of different policies regarding how service requests are handled, determining appropriate rules that will be applied to react to variable circumstances that the runway manager finds along the day.

This article focuses in the study of both services (landing and takeoff) that imply occupancy of the runway system of a single airport, taking into account the dynamically evolving circumstances found on a daily basis, represented by the original flight time schedule, and by the operational conditions that will be assumed by the runway system (airport runway configurations). We will study the effects of the different tactical decisions that are adopted throughout the day.

\section{Capacity and Operational Throughput of a runway system}
\begin{figure}
	\pgfkeys{/tikz/.cd,
		fa/.store in=\fa, 
		ru/.store in=\ru, 
		ex/.store in=\ex, 
		tt/.store in=\tt 
	}
	\begin{tikzpicture}[fa=5, ru=2, ex=1.5, tt=20, xscale=0.6, yscale=1.0]
		\draw[line width=0.5pt, -] (-0.5,{-\fa}) -- ({\tt},{-\fa}) node[above left] {approach};
		\draw[line width=0.5pt, -] (-0.5,0) -- ({\tt},0) node[above left] {threshold};
		\draw[line width=0.5pt, -] (-0.5,{\ex}) -- ({\tt},{\ex}) node[above left] {runway exit};
		\draw[line width=0.5pt, -] (-0.5,{\ru}) -- ({\tt},{\ru}) node[above left]{runway end};
		\draw[line width=6pt, -] (0.1,{\ex}) -- ({0.1+\ex-\ru},{\ru});
		\draw[line width=10pt, -] (0.1,0) -- (0.1,{\ru});
		\draw[white, dashed,line width=2pt, -] (0.1,0) -- (0.1,{\ru});
		\draw[white, dashed, line width=2pt, -] (0.1,{\ex}) -- ({0.1+\ex-\ru},{\ru});
		
		\def\t{2.5} 
		\def\s{2} 
		\draw[color=green,line width=1pt, variable=\xx,domain={0:2*\ex/\s},samples=50]
		plot ({\xx+\t},{\s*\xx*(4*\ex-\s*\xx)/(4*\ex)});
		\draw[color=green,line width=1pt]
		({\t},0) -- ({\t-\fa/\s},{-\fa});
		\draw[dotted,line width=1pt]
		({2*\ex/\s+\t},{\ex}) -- ({2*\ex/\s+\t},{0});
		\draw[gray, line width=4pt, -] ({\t},{0}) -- ({\t},{-3}) node[above,rotate=90]{safe};

		\def\t{2*\ex/2+2.5} 
		\def\s{1.8} 
		\draw[color=red,line width=1pt, variable=\xx,domain={0:2*\ru/\s+0.3},samples=50]
		plot ({\xx+\t},{(\s*\xx)^2/(4*\ru)});
		\draw[dotted,line width=1pt]
		({2*\ru/\s+\t},{\ru}) -- ({2*\ru/\s+\t},{0});
		\draw[gray, line width=4pt, -] ({\t},{0}) -- ({\t},{-2}) node[above,rotate=90]{safe};

		\def\t{6.5} 
		\def\s{1.3} 
		\draw[color=green,line width=1pt, variable=\xx,domain={0:2*\ex/\s},samples=50]
		plot ({\xx+\t},{\s*\xx*(4*\ex-\s*\xx)/(4*\ex)});
		\draw[color=green,line width=1pt]
		({\t},0) -- ({\t-\fa/\s},{-\fa});
		\draw[dotted,line width=1pt]
		({2*\ex/\s+\t},{\ex}) -- ({2*\ex/\s+\t},{0});
		\draw[gray, line width=4pt, -] ({\t-\fa/\s},{-\fa}) -- ({\t-\fa/\s},{3-\fa}) node[above,rotate=90]{safe};
		\draw[gray, line width=4pt, -] ({\t},{0}) -- ({\t},{-3}) node[above,rotate=90]{safe};
		
		\def\t{8.8} 
		\def\s{1.6} 
		\draw[color=green,line width=1pt, variable=\xx,domain={0:2*\ex/\s},samples=50]
		plot ({\xx+\t},{\s*\xx*(4*\ex-\s*\xx)/(4*\ex)});
		\draw[color=green,line width=1pt]
		({\t},0) -- ({\t-\fa/\s},{-\fa});
		\draw[dotted,line width=1pt]
		({2*\ex/\s+\t},{\ex}) -- ({2*\ex/\s+\t},{0});
		\draw[gray, line width=4pt, -] ({\t-\fa/\s},{-\fa}) -- ({\t-\fa/\s},{3-\fa}) node[above,rotate=90]{safe};
		\draw[gray, line width=4pt, -] ({\t},{0}) -- ({\t},{-3}) node[above,rotate=90]{safe};

		\def\t{2*\ex/1.6+8.8} 
		\def\s{2.2} 
		\draw[color=red,line width=1pt, variable=\xx,domain={0:2*\ru/\s+0.3},samples=50]
		plot ({\xx+\t},{(\s*\xx)^2/(4*\ru)});
		\draw[dotted,line width=1pt]
		({2*\ru/\s+\t},{\ru}) -- ({2*\ru/\s+\t},{0});
		\draw[gray, line width=4pt, -] ({\t},{0}) -- ({\t},{-2}) node[above,rotate=90]{safe};

		\def\t{2*\ex/1.6+2*\ru/2.2+8.8} 
		\def\s{1.6} 
		\draw[color=red,line width=1pt, variable=\xx,domain={0:2*\ru/\s+0.3},samples=50]
		plot ({\xx+\t},{(\s*\xx)^2/(4*\ru)});
		\draw[dotted,line width=1pt]
		({2*\ru/\s+\t},{\ru}) -- ({2*\ru/\s+\t},{0});
		\draw[gray, line width=4pt, -] ({\t},{0}) -- ({\t},{-2}) 	node[above,rotate=90]{safe};

		\def\t{15.1} 
		\def\s{1.5} 
		\draw[color=green,line width=1pt, variable=\xx,domain={0:2*\ex/\s},samples=50]
		plot ({\xx+\t},{\s*\xx*(4*\ex-\s*\xx)/(4*\ex)});
		\draw[color=green,line width=1pt]
		({\t},0) -- ({\t-\fa/\s},{-\fa});
		\draw[dotted,line width=1pt]
		({2*\ex/\s+\t},{\ex}) -- ({2*\ex/\s+\t},{0});
		\draw[gray, line width=4pt, -] ({\t-\fa/\s},{-\fa}) -- ({\t-\fa/\s},{3-\fa}) node[above,rotate=90]{safe};
	\end{tikzpicture}		
\caption{Space-time diagram of evolution for several flights using the runway according to safety distance restrictions (adapted from \cite{Newell_AirportCapacityandDelays}).}	\label{timespacediagram}
\end{figure}
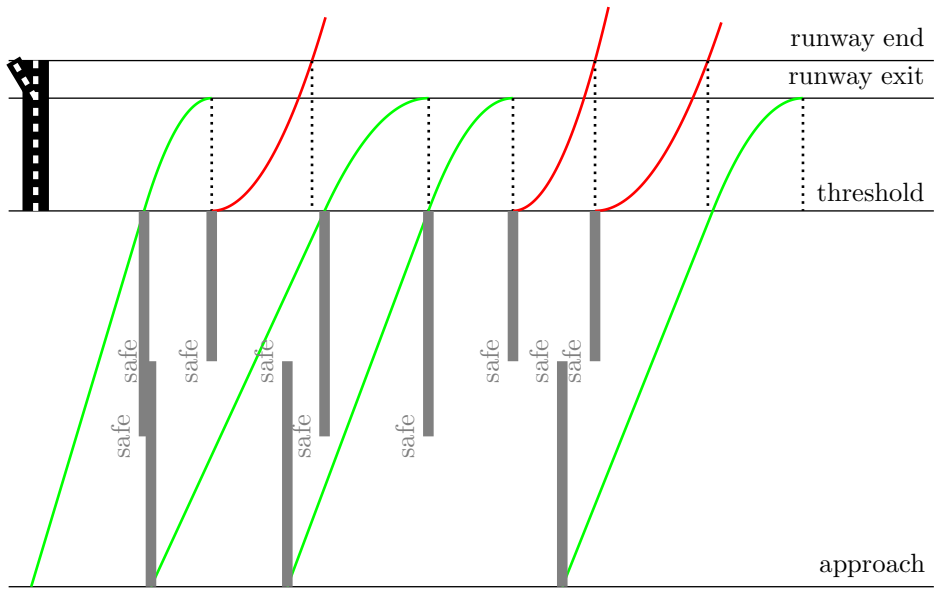

In \cite{Newell_AirportCapacityandDelays} Newell presents  a critical review of (pre-1980's) articles dealing with analytical models of runway operations. For a single runway, under stable conditions, all operations (landing and take-off) are performed along an imaginary line, always in the same (upwind) direction. There is a common path on flight (final approach) to be used by all airplanes for landing, a runway stretch on ground (from runway threshold to runway exit) to be used by all operations (landing and takeoff), and finally the takeoff path, which includes a final stretch of the runway and a neighboring airspace to be used only by takeoff airplanes. The evolution of several flights using the runway, by airplanes with different speeds, can be represented on a time-space diagram (figure \ref{timespacediagram}). Security concerns demand that at any moment the separation between airplanes is restricted: there is a prescribed separation between consecutive planes when entering or leaving the final approach sector; another separation for any airplane taking off, with respect to those in the final approach, and at any time the runway can be occupied by a single airplane. 

The variety of service times depending on the airplane, the possibility to start any landing or takeoff operation as long as there is no interference on the runway or at the common approach path, and other aspects related to these operations (visibility, random incidences, etc) make it possible to choose from diverse landing/takeoff sequences without infringement of security rules. Hence there is a possibility to choose a bidimensional parameter (number of landing and takeoff services) in a given time period (slot) so that the values of this bidimensional parameter lie within the operational capabilities of the runway system. It becomes relevant to fix a simplified model that handles the complexity of all factors and renders a sufficiently meaningful notion of capacity for an airport without need to involve an excessive number of parameters that, in general, would have a random nature making them hard to control.

For a modern perspective on the nature of landing and takeoff services, for the most relevant parameters used for its characterization and for a review of modern literature on the subject, one may consult \cite{Airport_runway_Scheduling_Bennet,rsp_survey_author_version}, where one finds different mathematical techniques that model these operations, together with several algorithms used in its optimization.

A first practical inference that one gets observing space-time diagrams (figure \ref{timespacediagram}) is that, though for a single class of operations the runway occupancy time depends linearly on the number of operations, when we alternate different classes of operations, the resulting runway occupancy time does not show this linear dependence. Based on these considerations Newell \cite{Newell_AirportCapacityandDelays} performed an analytical study determining the possibility to intercalate different classes of operations and airplanes so that the total number of operations can be tuned by a corresponding adjustment of the ratio between takeoff and landing airplanes. This analytical consideration leads to the result that the number $n^a,n^d\geq 0$ of arrivals (landings) and departures (takeoffs) that can be achieved in a certain time interval are restricted by a set of $J$ linear constraints of the type \begin{equation}
\label{Newellconstraints}	\alpha^a_j\cdot n^a_j+\alpha^d_j\cdot n^d_j\leq \beta_j \qquad (j=1\ldots J)
\end{equation} which together determine a convex domain that characterizes the throughput of the runway. This domain differs from airport to airport and may also differ taking into account specific circumstances, in particular those that impose particular operation modes and security rules due to meteorological conditions. Moreover, such a description is also possible not just for a single runway, but for the runway system of a given airport.

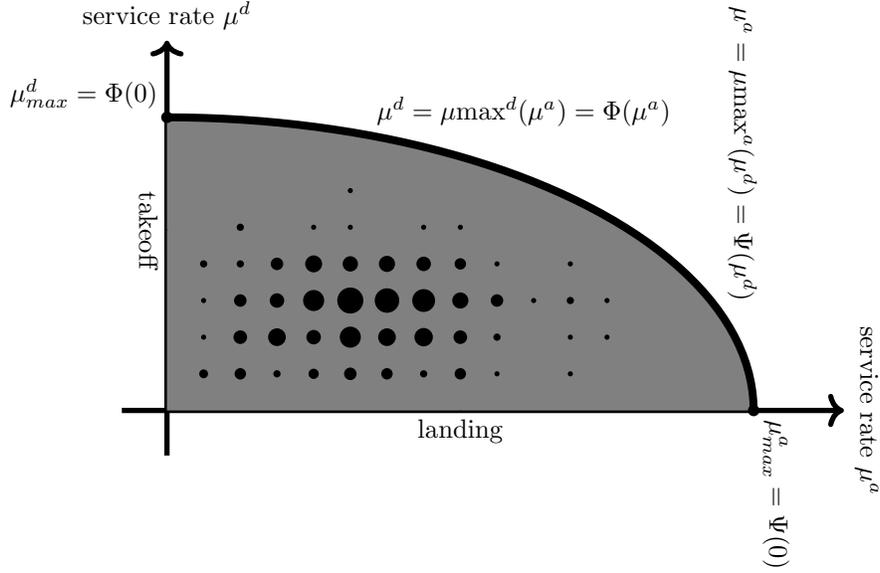
\begin{figure}
	\pgfkeys{/tikz/.cd,
	lamba/.store in=\lamba,
	lambd/.store in=\lambd,
	qa/.store in=\qa,
	qd/.store in=\qd,
}
\begin{tikzpicture}[lamba=1.9,lambd=1.2, qa=1.8, qd=2.9, scale=1.2, declare function = {
		z(\la,\mu,\q) = 1/\mu + (\q*\la)/(2*\mu*(\mu-\la));
		la(\mu,\z,\q) = \mu*2*(\z*\mu-1)/(\q+2*(\z*\mu-1));
		mu(\z,\la,\q)=0.5*(1+\la*\z+(1+\la^2*\z^2+2*\la*\z*(\q-1))^0.5)/\z;
	}
	]
	\draw[line width=2pt, ->] (-0.5,0) -- (7.5,0) node[above,rotate=-90] {service rate $\mu^a$};
	\draw[line width=2pt, ->] (0,-0.5) -- (0,4.1) node[above] {service rate $\mu^d$};
	
	
	\fill [gray, opacity=0.5, domain=0:90, variable=\t]
	(0,0)
	-- plot ({6.5*cos(\t)},{3.25*sin(\t)})
	-- cycle;

	\draw[line width=3pt,variable=\t,domain=0:90,samples=50]
	plot ({6.5*cos(\t)},{3.25*sin(\t))});
	\fill ({6.5*cos(70)},{3.25*sin(70)}) circle (0.0cm) node[above right] {$\mu^d=\mu\mathrm{max}^d(\mu^a)=\Phi(\mu^a)$};
	\fill ({6.5*cos(20)},{3.25*sin(20)}) circle (0.0cm) node[above left,rotate=-90] {$\mu^a=\mu\mathrm{max}^a(\mu^d)=\Psi(\mu^d)$};

	\fill (3.25,0) circle (0.0cm) node[below] {landing};
	\fill (0,2.) circle (0.0cm) node[below,rotate=-90] {takeoff};

	 \foreach \a / \b / \n in { 
	 	0.40625/	0.40625/	0.06,
	 	0.40625/	0.8125/	0.02,
	 	0.40625/	1.2188/	0.02,
	 	0.40625/	1.625/	0.04,
	 	0.8125/	0.40625/	0.1,
	 	0.8125/	0.8125/	0.14,
	 	0.8125/	1.2188/	0.12,
	 	0.8125/	1.625/	0.04,
	 	0.8125/	2.0312/	0.04,
	 	1.2188/	0.40625/	0.04,
	 	1.2188/	0.8125/	0.24,
	 	1.2188/	1.2188/	0.16,
	 	1.2188/	1.625/	0.12,
	 	1.625/	0.40625/	0.08,
	 	1.625/	0.8125/	0.16,
	 	1.625/	1.2188/	0.34,
	 	1.625/	1.625/	0.22,
	 	1.625/	2.0312/	0.02,
	 	2.0312/	0.40625/	0.12,
	 	2.0312/	0.8125/	0.34,
	 	2.0312/	1.2188/	0.52,
	 	2.0312/	1.625/	0.18,
	 	2.0312/	2.0312/	0.02,
	 	2.0312/	2.4375/	0.02,
	 	2.4375/	0.40625/	0.1,
	 	2.4375/	0.8125/	0.24,
	 	2.4375/	1.2188/	0.46,
	 	2.4375/	1.625/	0.2,
	 	2.8438/	0.40625/	0.04,
	 	2.8438/	0.8125/	0.26,
	 	2.8438/	1.2188/	0.4,
	 	2.8438/	1.625/	0.16,
	 	2.8438/	2.0312/	0.02,
	 	3.25/	0.40625/	0.1,
	 	3.25/	0.8125/	0.14,
	 	3.25/	1.2188/	0.2,
	 	3.25/	1.625/	0.1,
	 	3.25/	2.0312/	0.02,
	 	3.6562/	0.40625/	0.02,
	 	3.6562/	0.8125/	0.04,
	 	3.6562/	1.2188/	0.12,
	 	3.6562/	1.625/	0.02,
	 	4.0625/	1.2188/	0.02,
	 	4.4688/	0.40625/	0.02,
	 	4.4688/	0.8125/	0.02,
	 	4.4688/	1.2188/	0.04,
	 	4.4688/	1.625/	0.02,
	 	4.875/	0.8125/	0.02,
	 	4.875/	1.2188/	0.02}
	{
	\fill[color=black] ({\a},{\b}) circle [radius={0.2*\n^.5}];
	}	 

%
	
	
	\fill (0,3.25) circle (0.064cm) node[above left]{$\mu^d_{max}=\Phi(0)$};
	\fill (6.5,0) circle (0.064cm) node[above right,rotate=-90]{$\mu^a_{max}=\Psi(0)$};

\end{tikzpicture}
\caption{Operational throughput domain and envelope, and scatterplot of simulated operational throughput (numbers $n^a,n^d$ of arrivals and departures, circle area proportional to frequency) for different 15min time slots at a given airport with an specific configuration of its runway system.}\label{scatterplot}\end{figure}

Instead of making theoretical considerations about the runway system and the different configurations under which this system may operate, one may also look at historical data of the number of landing and takeoff operations that a given airport has achieved in several time slots, under certain meteorological conditions and runway configurations. This information $(n^a,n^d)$ is called operational throughput of the runway system in different time slots and can be represented by a 2D scatter plot (figure \ref{scatterplot}). 

There are different ways to define and measure an airport capacity. The origin of this research can be found in \cite{Blumstein_AnaliticalInvestigationAirportCapacity}. Using a loose description, we may say that a capacity measure should represent the number of movements (in its bidimensional nature, discriminating landings from takeoffs) that can be performed in a given time interval, when the runway system is taken to its limits without violation of security rules, and in the presence of a steady demand of service.

It was Gilbo \cite{Airport_Capacity_Gilbo93} the first in considering the compilation of historical landing and takeoff registers from an airport to infer the shape of the operational throughput envelope proposed by Newell. This led to the consideration of a domain on the positive quadrant of the $(n^a,n^d)$ plane, as the object that characterizes the capacity of a runway system operating under certain configuration. Gilbo used this domain as cornerstone in the management of a runway system seeking to maximize the performance on a single day, by choosing among the available ratios of landing vs takeoff. 

The determination of a convex capacity domain from historical data of a specific airport can be found in \cite{Airport_Capacity_Gilbo93} and more recent works \cite{SimaiakisThesis13_Gilbocapacities}. Using the operational throughput scatterplot showing the number of landings and takeoffs in several 15min slots, Gilbo studied the operational throughput envelope, a boundary of this domain determined by the maximal number of operations (after elimination of outliers) of one of these services, at the time-slots where we observe a given number of operations of the second type.

For a given configuration of runways, we may choose a service policy (or simply a policy), that is, service rates $(\mu^a,\mu^d)$ measured in (possibly non-integer) number of operations per slot. The application of this policy determines how operations should be performed  for a given time interval. Only certain policies are within the runway system capacities, in particular the service rates should lie on the positive quadrant $\mu^a\geq 0$, $\mu^d\geq 0$, and below a convex region limited by a curve, the ``operational throughput envelope''.

For a specific intended landing service rate $\mu^a$, there exist a maximal takeoff service rate $\mu\mathrm{max}^d(\mu^a)$ that could be achieved for the takeoff operations, while ensuring the rate $\mu^a$ for arrivals. These values  determine a plane curve in the positive quadrant $(\mu^a,\mu^d)$:
$$\mu^d=\mu\mathrm{max}^d(\mu^a)=\Phi(\mu^a)$$
This curve is called operational throughput envelope, and was described by Gilbo using the monotone decreasing concave function $\Phi$. Admissible policies $(\mu^a,\mu^d)$ are characterized by the restriction $\mu^d\leq \Phi(\mu^a)$ on the positive quadrant ($\mu^a,\mu^d\geq 0$). 

In the same way, for each specific takeoff throughput rate $\mu^d$ we may observe the maximal landing throughput rate $\mu\mathrm{max}^a(\mu^d)$. The operational throughput envelope is characterized by :
$$\mu^a=\mu\mathrm{max}^a(\mu^d)=\Psi(\mu^d)$$
where $\Psi$ is again a monotone decreasing concave function, inverse of the function  $\Phi$, that is: $\Psi=\Phi^{-1}$.

All operational throughput observations available in the historical registers for a given runway system that operates following a specific configuration are the result of applying an admissible policy in specific circunstances and should lie within the operational throughput domain (figure \ref{scatterplot}). 
Relevant values associated to this envelope are:
$$\mu\mathrm{max}^a(0)=\mu^a_{max},\qquad \mu\mathrm{max}^d(0)=\mu^d_{max}$$
which represent the maximal number of landings or takeoffs that can be performed with the given runway system configuration, and that would be achieved in optimal circumstances under a maximal priority policy of one of these services.

\begin{figure}
	\pgfkeys{/tikz/.cd,
		lamba/.store in=\lamba,
		lambd/.store in=\lambd,
		qa/.store in=\qa,
		qd/.store in=\qd,
	}
	\begin{tikzpicture}[lamba=1.9,lambd=1.2, qa=1.8, qd=2.9, scale=1.2, declare function = {
			z(\la,\mu,\q) = 1/\mu + (\q*\la)/(2*\mu*(\mu-\la));
			la(\mu,\z,\q) = \mu*2*(\z*\mu-1)/(\q+2*(\z*\mu-1));
			mu(\z,\la,\q)=0.5*(1+\la*\z+(1+\la^2*\z^2+2*\la*\z*(\q-1))^0.5)/\z;
		}
		]
		\draw[line width=2pt, ->] (-0.5,0) -- (7.5,0) node[above,rotate=-90] {service rate $\mu^a$};
		\draw[line width=2pt, ->] (0,-0.5) -- (0,4.1) node[above] {service rate $\mu^d$};
		
		
		\fill [gray, opacity=0.5, domain=0:90, variable=\t]
		(0,0)
		-- plot ({6.5*cos(\t)},{3.25*sin(\t)})
		-- cycle;

		\draw[line width=3pt,variable=\t,domain=0:90,samples=50]
		plot ({6.5*cos(\t)},{3.25*sin(\t))});
		\fill ({6.5*cos(70)},{3.25*sin(70)}) circle (0.0cm);
		\fill ({6.5*cos(20)},{3.25*sin(20)}) circle (0.0cm);
		
		\foreach \t/\j in {0/0,15/1,30/2,45/3,60/4,75/5,90/J} {
			\fill[color=red] ({6.5*cos(\t)},{3.25*sin(\t))}) circle [radius=0.1] node[above right] {$(x_{\j},y_{\j})$};
		}

		\draw[dashed,color=red,line width=2pt] ({6.5*cos(15)},0) \foreach \t in {0,15,30,45,60,75,90} {
			--({6.5*cos(\t)},{3.25*sin(\t)})} -- (0,{3.25*sin(90)});
		
		\fill (0,3.25) circle (0.064cm) node[above left]{$\mu^d_{max}=\Phi(0)$};
		\fill (6.5,0) circle (0.064cm) node[above right,rotate=-90]{$\mu^a_{max}=\Psi(0)$};
		
			\fill (3.25,0) circle (0.0cm) node[below] {landing};
		\fill (0,2.) circle (0.0cm) node[below,rotate=-90] {takeoff};
		
	\end{tikzpicture}
	\caption{Operational throughput envelope and its approximation as convex polygon  using the sequence of control points (\ref{ordercontrolpoints}).}	\label{controlpoints}
\end{figure}
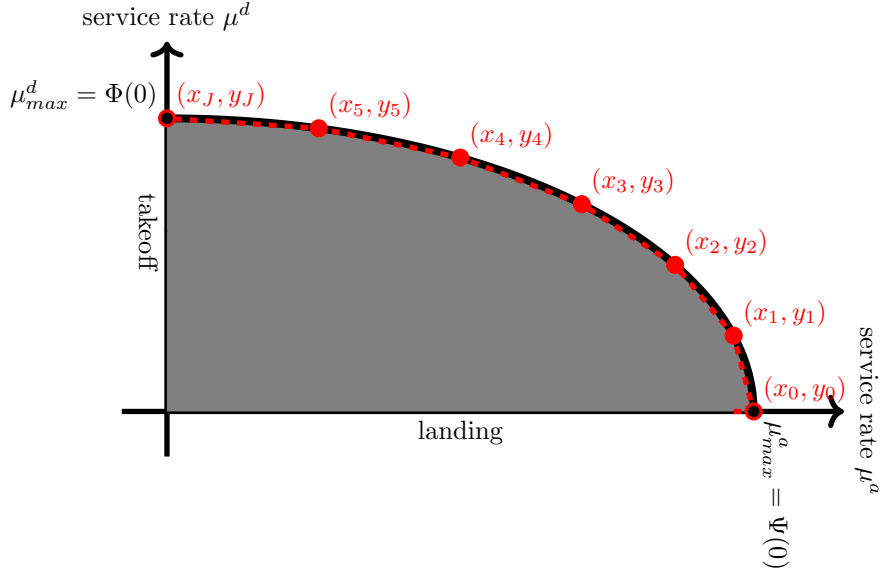

Such domains have been identified for a large number of airports (see the collection for major US airports available in \cite{USenvelopes}). Their boundary is usually given by both axis and a polygonal line limiting a convex domain for $(\mu^a,\mu^d)$. 	The vertices of this convex domain are $(0,0)$ together with a sequence of points (control points) ordered as follows:
	\begin{equation}\label{ordercontrolpoints}
	\begin{aligned}
	\left( (x_j,y_j)\right)_{j=0\ldots J},\text{ decreasing }(x_j),\text{ increasing } (y_j),\, x_J=0,\, y_0=0
	\end{aligned}
	\end{equation}
as shown in figure \ref{controlpoints}. This usually represents the only known information $\Phi(x_j)=y_j$, $\Psi(y_j)=x_j$ for Gilbo's envelope function $\Phi$, which is empirical and has no specific analytical form.

 Using the vertices of the polygonal line we get a description of Gilbo's domain for policies $(\mu^a,\mu^d)$ similar to the form (\ref{Newellconstraints}) that constrained effective operational throughput $(n^a,n^d)$ :
	\begin{equation}\label{segments}\begin{aligned}
	&\mu^a\geq 0,\quad \mu^d\geq 0\\ 
	\mu^a\cdot (y_j-y_{j-1})+\mu^d\cdot (x_{j-1}-x_j) &\leq x_{j-1}\cdot y_j-x_j\cdot y_{j-1},\quad j=1\ldots J
	\end{aligned}
	\end{equation}
We will assume that the capacity of a runway system under a specific configuration is described by its associated operational throughput envelope, and that this envelope is described as a convex polygon with vertices at $(0,0)$ and an ordered sequence of control points (\ref{ordercontrolpoints}).

Applying a policy $(\mu^a,\mu^d)$ for a given time slot does not imply that the expected number of landing or takeoff operations will coincide with the rates $(\mu^a,\mu^d)$ declared as policy. As we said, these values are somehow maximal values for the number of operations that, historically, have been observed under particular (normally exceptional) circumstances. The number of operations that are effectively performed depends in particular on the number of clients (airplanes) demanding those services, at the given time slot, a situation that is best described from the perspective of queuing theory.

\section{Runway operation as competing queue systems}\label{sec3}

As observed by Shone et al. \cite{Queue_theory_runway_published.pdf,OR60_paper}, airport runway operations is a specific instance of the problem of distribution of a common resource between two sources of clients that have a variable demand along time. In this case the resource is the runway occupancy time and sources are the set of airplanes trying to land or takeoff. The level of demand is determined by the service schedule each day. This situation may be studied within the general model of queue theory. Gilbo (1993) \cite{Airport_Capacity_Gilbo93},  Bertsimas Frankowich and Odoni (2011) \cite{BertFranOdon11, Frankowich12_TFModelAirports}, del-Olmo, Lulli \cite{Olmo_Lulli} or Gluschenko \cite{Dynamic_usage_of_capacity_for_arrivals_and_departures_in_queue_minimization} started from a deterministic demand rate, with scheduled client arrivals that where not affected by random factors.  
Jacquillat and Odoni \cite{JacOdoni15GoodQueue} and  Ignacollo \cite{A_Simulation_model_for_airport_capacity_and_delay_analysis} on the other hand recovered Newell's model from a perspective of stochastic service procurement, closer to the reality. While Newell assumed that the service time follows a Poisson model,  Jacquillat, Odoni e Ignacollo assumed it to follow Erlang's model. Using these models is partially justified by the nature of the service but mostly for being models where the queue sizes can be described in a precise way. There are empirical studies on other possible choices for the probability distribution of takeoff service times, by Simaiakis and Balakrishnan \cite{SimaiakisBalakrishnanISIATM2013} or by Pujet et al. \cite{PujDelFer_departure_process.pdf}.

The difference between these two assumptions (deterministic or stochastic) is important. It is intuitively understood that if the needed time span to execute a service is larger than the expected time between client arrivals, the whole system will generate a queue that grows with time, without any bound on its size. Less intuitive is the situation where the time span to execute a service lies below the expected time between client arrivals: if the arrival is deterministic and uniform in time, this system will never generate a queue; however, if the inter-arrival time depends on random factors, on average the system will have a queue and the expected size of this queue stabilizes as time advances. In this case, clients arriving at the airport should count on a certain time needed for the operation (which was considered by Gilbo, Bertsimas or Frankowich), plus an additional waiting time that in fact depends on the demand observed for that service. As stressed by Kim and Hansen \cite{Kim} the relationship between throughput and demand at airports, is only roughly approximated by a model based on the conventional notion of capacity

Basic queuing models state that a system with an infinite source of clients that demand a service at random times and a single server that delivers the service by order of arrival, can be characterized by a sequence of positive random variables $C_k$ measuring the time between the arrival of the $k$th and $(k+1)$th client, and positive variables $S_k$ measuring the time expended to serve the $k$th client. It is natural to assume that all variables $C_k$ correspond to the same continuous probabilistic model with finite expected value and coefficient of variation, hence represented as a variable $C$ with expected value $\EE(C)=1/\lambda>0$ and with finite variance $$\mathop{\mathrm{Var}}(C)=\EE\left((C-1/\lambda)^2\right)=\EE(C^2)-1/\lambda^2=\frac{q_C-1}{\lambda^2}$$ where we denote $q_C=\EE(C^2)/(\EE(C))^2\geq 1$ the quadratic ratio of momenta, related to the coefficient of variation $c$ by $q_C=1+c^2$. It is also natural to assume that all variables $S_k$ are represented by a variable $S$ from a continuous probabilistic model with expected value $\EE(S)=1/\mu$ and quadratic ratio of momenta $\EE(S^2)/(\EE(S))^2=q_S\geq 1$. Due to the dimensionless nature of the quadratic ratios of momenta, we may assume that a given server will present a characteristic value of $q_S$, which doesn't change even in cases where the service rate is variable.

Associated to such a system there is a natural stochastic process $Q(t)$, a time-parameterized family of positive integer random variables that measure, at each instant $t$, the queue length of clients waiting for the service. The laws governing this process are described in terms of both variables $C$ and $S$. A stationary system corresponds to the case in which there is a fixed random variable $Q$ with specific distribution such that $Q(t)=Q$ for each time $t$. The expected value $\EE(Q)$ is then called the expected queue length in the stationary situation. If $Q_0=Q(0)$ is not stationary, it will evolve and for $t\to \infty$ may converge (in probability) to some distribution $Q_{\infty}$ (we say that the system stabilizes), leading to an expected value $\EE(Q_{\infty})$, the long-term expected queue length.  

Fundamental parameters in the system are the client arrival rate $\lambda$ (inverse of the expected time between two consecutive arrivals), and the client service rate $\mu$ (inverse of the expected time that takes the service of a client). At any time the quotient $\rho=\frac{\lambda}{\mu}$ between arrival and service rates is called the utilization rate of the service. Both quadratic ratios of momenta $q_S,q_C$ are also relevant parameters (accounting for second order momenta of random variables $C,S$) associated to these variables.

\begin{figure}
	\pgfkeys{/tikz/.cd,
		mua/.store in=\mua,
		qa/.store in=\qa,
		pa/.store in=\pa
	}
	\begin{tikzpicture}[mua=10, qa=1.5, pa=0.333, xscale=0.9, yscale=0.9, declare function = {
			z(\la,\mu,\q) = 1/\mu + (\q*\la)/(2*\mu*(\mu-\la));
			la(\mu,\z,\q) = \mu*2*(\z*\mu-1)/(\q+2*(\z*\mu-1));
			mu(\z,\la,\q)=0.5*(1+\la*\z+(1+\la^2*\z^2+2*\la*\z*(\q-1))^0.5)/\z;
		}
		]
		\draw[line width=2pt, ->] (-0.5,0) -- (11,0) node[above, rotate=-90] {$\lambda$ (arrivals/slot)};
		\draw[line width=2pt, ->] (0,-1) -- (0,{2*z(9.7,\mua,\qa)}) node[above] {$z$ (slots of delay)};
		
		%
		
		\draw[line width=1pt,variable=\la,domain=0:9.7,samples=50]
		plot ({\la},{2*z(\la,\mua,\qa)});
		\fill ({0},{2*\pa}) circle (0.1cm) node[above right] {predefined $p$};
		\fill ({\mua},{0}) circle (0.1cm) node[below] {$\lambda_{sta}$};
		\fill ({la(\mua,\pa,\qa)},{0}) circle (0.1cm) node[below] {$\lambda_{cong}$};
		
		\fill ({la(\mua,\pa,\qa)/2},{0.1+2*\pa}) circle (0.cm) node[below left] {sustainable};
		\fill ({(la(\mua,\pa,\qa)+\mua)/2},{2}) circle (0.cm) node[rotate=60] {congestion};
		\fill (11.5,4) circle (0.cm) node[above,rotate=90] {saturation};
		
		\draw[dashed] (0,{2*\pa}) -- ({la(\mua,\pa,\qa)},{2*\pa})-- ({la(\mua,\pa,\qa)},0);

		\fill [green, opacity=0.5, domain={0:la(\mua,\pa,\qa))}, variable=\la]
		(0,0)-- (0,{2*\pa})
		-- plot ({\la},{2*z(\la,\mua,\qa)})
		-- ({la(\mua,\pa,\qa)},0) -- cycle;	
		
		\fill [yellow, opacity=0.5, domain={la(\mua,\pa,\qa)):9.7}, variable=\la]
		({la(\mua,\pa,\qa)},0)-- ({la(\mua,\pa,\qa)},{2*\pa})
		-- plot ({\la},{2*z(\la,\mua,\qa)})
		-- (\mua,{2*z(9.7,\mua,\qa)}) -- (\mua,0) -- cycle;	
		
		\fill [red, opacity=0.5]
		(\mua,0) -- (\mua,{2*z(9.7,\mua,\qa)}) -- (11,{2*z(9.7,\mua,\qa)}) -- (11,0) --cycle;	
		

		%
		
	\end{tikzpicture}
	\caption{For a given service rate $\mu$, certain arrival rates $\lambda$ may lead to stable queue systems with large transit times (congestion) or reasonable transit times (under a predefined level $p$), or may produce a non-controlled growth in transit times (saturation).}	\label{figsaturation} 
\end{figure}
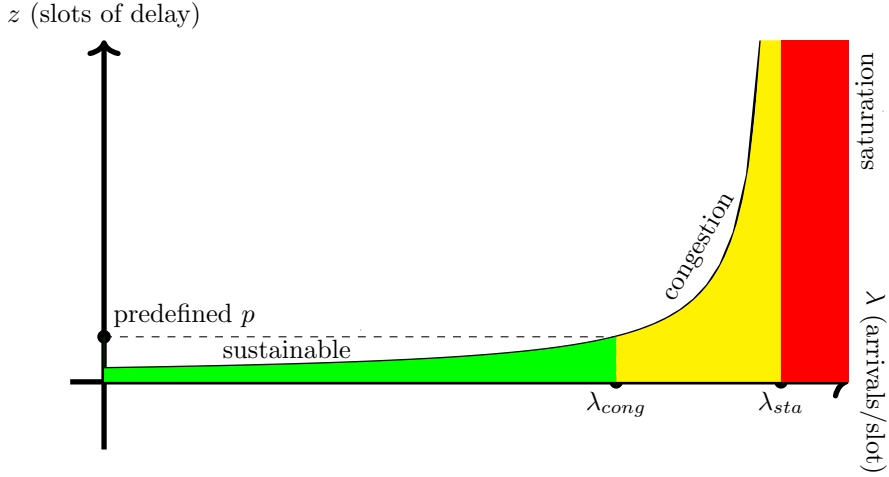

For utilization rate $\rho\geq 1$ no stationary solution exists (we say the system is unstable) and there exist ever growing expected values for the number of clients at the queues, in a situation that we call saturation conditions for the system. An utilization rate $\rho$ slightly below 1 allows for states that are stationary, with queue size distributions that have expected values with an (possibly large) uniform bound for all time. The large number of clients in the queue, however, makes us talk of congestion conditions for the system. Only utilization rates $\rho$ below a certain predefined level represent the buildup of small queues, which are stable, and shall be called sustainable conditions for the system. Queue size for saturated services tend to unlimited growth, and queue size for congestion conditions is very sensitive to small variations on the parameters and are considered as non-sustainable. Large queue sizes for runway services at any airport may severely influence all remaining services provided by the airport. Therefore, runway (landing or takeoff) services  are usually restricted to operate under a predefined congestion level (see figure \ref{figsaturation}).

The celebrated Pollakzec-Khintchine formula states that in the case where $C$ is a variable of the exponential model, the stationary solution $Q$ for the queue length has the following expected value:
$$\EE(Q)=\frac{\rho^2(c_S^2+1)}{2(1-\rho)}=\frac{\lambda^2\cdot q_S}{2\mu\cdot (\mu-\lambda)}$$
The general situation (for continuous random variables $C,S$ from arbitrary models) is more difficult to deal with, but an estimation \cite{Kingman} for the long-term expected value of the queue length is:
	\begin{equation}\label{Kingmanestimation}
		\EE(Q_{\infty})=\frac{\rho^2}{1-\rho}\cdot \frac{q_S+q_C-2}{2}=\frac{\lambda^2}{\mu\cdot(\mu-\lambda)}\cdot \frac{q_S+q_C-2}{2}
		\end{equation}
which coincides with Pollakzec-Khintchine formula in the case $q_C=2$ (recall that the exponential model has coefficient of variation 1 and thus a quadratic ratio of momenta 2)

This formula also leads to a value for waiting times using Little's law \cite{Littles-Law-50-Years-Later}:
$$\EE(W_{\infty})=\frac{1}{\lambda}\EE(Q_{\infty})=\frac1{\mu\rho} \cdot \frac{\rho^2}{1-\rho}\cdot \frac{q_S+q_C-2}{2}=\frac{\lambda}{\mu(\mu-\lambda)}\cdot \frac{q_S+q_C-2}{2}$$
If we add the expected time for a client to be served after waiting in the queue, we conclude for the transit time $Z=S+W_{\infty}$ that each client spends in the system:
\begin{equation}\label{calculaz}
	z=\EE(Z)=\frac1\mu\cdot \left(1+\frac{\rho}{1-\rho}\cdot\frac{q}{2}\right)=\frac1\mu\cdot \left(1+\frac{q\cdot \lambda}{2(\mu-\lambda)}\right)
\end{equation}
where the coefficient $q=q_S+q_C-2\geq 0$ is 0 for the deterministic arrival and service of clients, it is $q=1$ for Poisson arrival processes and deterministic service time, and it is $q\geq 1$ for Poisson arrival processes and random service times. The stable transit time $z$ measures the delay experienced by each client that arrives at a stationary state of the system, and is a natural measure of performance for the system.

For these queuing systems that are not unstable, there are constraints (\ref{calculaz}) that relate arrival and service rates $\lambda$, $\mu$ of clients (where $0<\lambda<\mu$) with the expected time of transit $z>\frac1\mu$ for each client that arrives in the system in the stable situation (the stable transit time). Solving (\ref{calculaz}) in the specific case $q>1$ for each of the variables in terms of the remaining ones:
	\begin{equation}\label{relatezmulambda}\begin{aligned} z&=\frac{1}{\mu}+\frac{q\lambda}{2\mu(\mu-\lambda)}\\ \mu&=\frac{1+\lambda z+\sqrt{1+\lambda^2z^2+2\lambda z(q-1)}}{2z}\\ \lambda&=\frac{2(z\mu-1)}{q+2(z\mu-1)}\cdot \mu\end{aligned}
	\end{equation}

For any given client service rate $\mu>0$, the arrival rates $\lambda\geq 0$ above a given value $\lambda_{sta}=\mu$ lead to saturation conditions for the system. If a predefined level of service $p$ is fixed (a delay tolerance), we also split stable situations (we mean non-saturated) into two cases,  those with arrival rates above a given level $\lambda_{cong}$ called congestion conditions, and those with arrival rates below this level, which we call sustainable conditions.

To avoid the congestion or saturation of the airport, there may be a predefined delay tolerance, fixing an upper bound for the queue size that we may expect, and this restriction can be described, in terms of delays associated to flights waiting to be served rather than of the number of these flights. The utilization rate that leads to such stable small (below predefined delay tolerance) queues can be maintained for long periods of time and can be used as a sustainable capacity of the airport, not to be mistaken with the maximal capacity, that is associated to saturation and congestion conditions. This situation is illustrated in figure \ref{figsaturation}.

Both landing and takeoff services for a runway system that operates with occupancy rates below 1 and in constant conditions have the previously defined behavior and, in the long term, we may expect the stabilization of the properties of queues that will appear for both types of operations. For simplicity, time will be measured in a standard unit, the slot size, and arrival and service of clients will be given in airplanes per slot.

\section{Sustainable policies}\label{sec4}

Assume a certain runway system configuration is used for a certain time slot. The runway system becomes a shared server for both landing and takeoff operations where clients are airplanes demanding either of these services. Assume the slot has  a schedule (subject to randomness) with $\lambda^a$ landing and $\lambda^d$ takeoff operations as expected demand (we will use indices a--``arrival'' or d--``departure''). We have two queue systems attending clients with characteristic variables $C^a,S^a,C^d,S^d$ that we assume to be independent. Parameters $(\lambda^a,\lambda^d)$ may be used as demand rates for the services of landing and takeoff, respectively, if we use the slot as time unit. 

Assume that a policy $(\mu^a,\mu^d)$ should be chosen for the service rate of arrivals and departures. The policy should lie on the operational throughput domain. 
Observe, however, that for a given landing and takeoff schedule $(\lambda^a,\lambda^d)$ certain policies $(\mu^a,\mu^d)$ within the operational throughput envelope may lead to either saturation or congestion conditions for the landing or takeoff services. Saturation happens whenever $\mu^a\leq \lambda^a$ or $\mu^d\leq \lambda^d$ . In the case of no saturation, the services will have a corresponding queue that in the long term renders stable transit times $z^a\geq 1/\mu^a$, $z^d\geq 1/\mu^d$ given by formulas (\ref{relatezmulambda}) and interpreted as flight delays.

Congestion (stable transit times above a predefined level) appears when $\mu^a<\mu^a_{cong}$ or $\mu^d<\mu^d_{cong}$, where these congestion levels $\mu^a_{cong}$, $\mu^d_{cong}$ depend on the predefined service levels $p^a,p^d$ following formulas (\ref{relatezmulambda}):
	\begin{equation}\label{mucong}
	\begin{aligned} \mu^a_{cong}&=\frac{1+\lambda^a p^a+\sqrt{1+(\lambda^a)^2(p^a)^2+2\lambda^a p^a(q^a-1)}}{2p^a}\\
	\mu^d_{cong}&=\frac{1+\lambda^d p^d+\sqrt{1+(\lambda^d)^2(p^d)^2+2\lambda^d p^d(q^d-1)}}{2p^d}
	\end{aligned}\end{equation}
Hence when a certain demand rate $(\lambda^a,\lambda^d)$ and a predefined service level $(p^a,p^d)$ are given, Gilbo's operational throughput domain and its envelope are divided into regions that identify policies that lead to saturation or to congestion conditions of either service, or to sustainability conditions for both services, as seen in figure \ref{figuraGilboDecomposto}.

\pgfkeys{/tikz/.cd,
	lamba/.store in=\lamba,
	lambd/.store in=\lambd,
	qa/.store in=\qa,
	qd/.store in=\qd,
	pd/.store in=\pd,
	pa/.store in=\pa
}
	\begin{figure}
	\begin{tikzpicture}[lamba=2.,lambd=1.2, qa=2.1, qd=4.2, pa=1.4, pd=2.7, scale=1.35, declare function = {
		z(\la,\mu,\q) = 1/\mu + (\q*\la)/(2*\mu*(\mu-\la));
		la(\mu,\z,\q) = \mu*2*(\z*\mu-1)/(\q+2*(\z*\mu-1));
		mu(\z,\la,\q)=0.5*(1+\la*\z+(1+\la^2*\z^2+2*\la*\z*(\q-1))^0.5)/\z;
		valory(\valorx)=0.5*(6.5^2-\valorx^2)^0.5;
		valorx(\valory)=(6.5^2-4*\valory^2)^0.5;
	}
	]
	\draw[line width=2pt, ->] (-0.5,0) -- (7.1,0) node[right] {$\mu^a$};
	\draw[line width=2pt, ->] (0,-0.5) -- (0,4.1) node[left] {$\mu^d$};
	\fill (0,3.25) circle (0.064cm) node[above left]{$\mu^d_{max}=\Phi(0)$};
	\fill (6.5,0) circle (0.064cm) node[above right,rotate=-90]{$\mu^a_{max}=\Psi(0)$};
	\draw[dashed] (0,{\lambd}) node[left]{$\lambda^{d}$}-- ({valorx(\lambd)},{\lambd})-- ({valorx(\lambd)},0);
	\draw[dashed] (0,{mu(\pd,\lambd,\qd)}) node[left]{$\mu^{d}_{cong}$}-- ({valorx(mu(\pd,\lambd,\qd))},{mu(\pd,\lambd,\qd)})-- ({valorx(mu(\pd,\lambd,\qd))},0);
	
	\draw[dashed] ({\lamba},0) node[below]{$\lambda^{a}$}-- ({\lamba},{valory(\lamba)})--(0,{valory(\lamba)});
	\draw[dashed] ({mu(\pa,\lamba,\qa)},0) node[right, rotate=-90]{$\mu^a_{cong}$}-- ({mu(\pa,\lamba,\qa)},{valory(mu(\pa,\lamba,\qa))})--(0,{valory(mu(\pa,\lamba,\qa))});
	
	
	\fill [red, opacity=0.5, domain=0:{asin(\lambd/3.25)}, variable=\t]
	({\lamba},0)-- plot ({6.5*cos(\t)},{3.25*sin(\t)})-- ({\lamba,\lambd})--cycle;	
	\fill [red, opacity=0.5, domain={acos(\lamba/6.5)}:90, variable=\t]
	(0,0)-- ({\lamba},0)-- plot ({6.5*cos(\t)},{3.25*sin(\t)})-- cycle;	
	
	\fill [yellow, opacity=0.5, domain={acos(\lamba/6.5):asin(\lambd/3.25)}, variable=\t]
	({\lamba}, {\lambd})
	-- plot ({6.5*cos(\t)},{3.25*sin(\t)})
	-- cycle;	
	
	\fill [green, opacity=0.5, domain={acos((mu(\pa,\lamba,\qa))/6.5):asin((mu(\pd,\lambd,\qd))/3.25)}, variable=\t]
	({mu(\pa,\lamba,\qa)}, {mu(\pd,\lambd,\qd)})
	-- plot ({6.5*cos(\t)},{3.25*sin(\t)})
	-- cycle;

	\draw[line width=3pt, variable=\t,domain={0:90},samples=50]	plot ({6.5*cos(\t)},{3.25*sin(\t))});

	\fill ({\lamba},{valory(\lamba)}) circle (0.0cm);
	\fill ({valorx(\lambd)},{\lambd}) circle (0.0cm);
	\foreach \t in {15,30,45,60,75} {
		\fill[color=black] ({6.5*cos(\t)},{3.25*sin(\t))}) circle [radius=4pt];
	}
	\foreach \t in {0,15,75,90} {
		\fill[color=red] ({6.5*cos(\t)},{3.25*sin(\t))}) circle [radius=3pt];
	}
	\foreach \t in {45,60} {
		\fill[color=green] ({6.5*cos(\t)},{3.25*sin(\t))}) circle [radius=3pt];
	}
	\fill[color=yellow] ({6.5*cos(30)},{3.25*sin(30))}) circle [radius=3pt];
	\fill[color=blue]{(\lamba,\lambd)} circle [radius=4pt];
	\draw[color=red,line width=2pt,variable=\t,domain={acos(\lamba/6.5)}:90,samples=50]
	plot ({6.5*cos(\t)},{3.25*sin(\t)});
	\draw[color=red,line width=2pt,variable=\t,domain=0:{asin(\lambd/3.25)},samples=50]
	plot ({6.5*cos(\t)},{3.25*sin(\t)});
	\draw[color=yellow,line width=2pt, variable=\t, domain={acos((mu(\pa,\lamba,\qa))/6.5)}:{acos(\lamba/6.5)},samples=50]
	plot ({6.5*cos(\t)},{3.25*sin(\t)});
	\draw[color=yellow,line width=2pt, variable=\t, domain={asin(\lambd/3.25)}:{asin((mu(\pd,\lambd,\qd))/3.25)},samples=50]
	plot ({6.5*cos(\t)},{3.25*sin(\t)});

	\draw[color=green,line width=2pt, variable=\t, domain={acos((mu(\pa,\lamba,\qa))/6.5)}:{asin((mu(\pd,\lambd,\qd))/3.25)},samples=50]
	plot ({6.5*cos(\t)},{3.25*sin(\t)});
	
	
	
	\fill ({valorx(\lambd)},0) circle (0.064cm) node[right, rotate=-90]{$\mu^a_{sta}=\Psi(\lambda^d)$};
	\fill ({valorx(mu(\pd,\lambd,\qd))},0) circle (0.064cm) node[below right,rotate=-90]{$\Psi(\mu^d_{cong})$};

	\fill (0,{valory(\lamba)}) circle (0.064cm) node[left] {$\mu^d_{sta}=\Phi(\lambda^a)$};
	\fill (0,{valory(mu(\pa,\lamba,\qa))}) circle (0.064cm) node[below left] {$\Phi(\mu^a_{cong})$};
	
	\draw [black] plot [only marks, mark size=2.5, mark=square*] coordinates {({mu(\pa,\lamba,\qa)},{valory(mu(\pa,\lamba,\qa))})} node[above right] {$(\bar x_J,\bar y_J)$};

	\draw [black] plot [only marks, mark size=2.5, mark=square*] coordinates { ({valorx(mu(\pd,\lambd,\qd))},{mu(\pd,\lambd,\qd)}) } node[above right] {$(\bar x_0,\bar y_0)$};

	\draw [black] plot [only marks, mark size=2.5, mark=square*] coordinates { ({mu(\pa,\lamba,\qa)},{mu(\pd,\lambd,\qd)})};
	
\end{tikzpicture}
	\caption{For a given demand rate and a predefined delay tolerance, Gilbo's domain splits in different regions, where the situation is sustainable, congestion or saturation, for the pair of landing/takeoff services.}\label{figuraGilboDecomposto}
\end{figure}
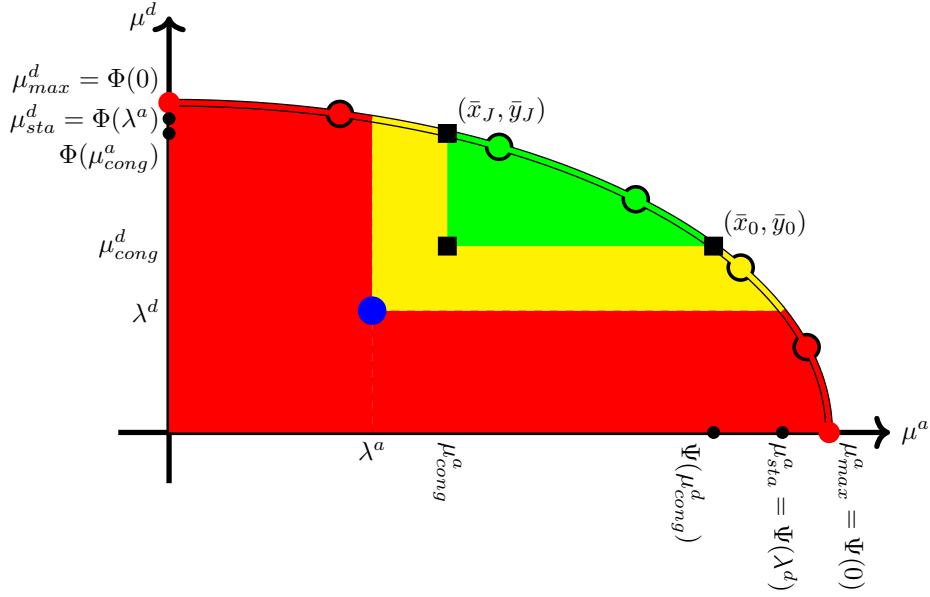

To avoid saturation for the landing service, the policy must have a service rate $\mu^a$ for landing above the corresponding demand rate, $\mu^a>\lambda^a$, which implies that the service rate $\mu^d$ for takeoff should remain below the rate $\mu^d_{sta}=\Phi(\lambda^a)$.

To avoid saturation for the takeoff service, the policy must have a service rate $\mu^d$ for takeoff above the corresponding demand rate, $\mu^d>\lambda^d$, which implies that the service rate $\mu^a$ for landing should remain below the rate $\mu^a_{sta}=\Psi(\lambda^d)$.
 
In a similar way, to avoid congestion for the landing service, the policy must have at least $\mu^a_{cong}$ as service rate for landing, $\mu^a\geq \mu^a_{cong}$, which implies at most $\Phi(\mu^a_{cong})$ as service rate for takeoff, that is, $\mu^d\leq \Phi(\mu^a_{cong})$.

To avoid congestion for the takeoff service, the policy must have a service rate $\mu^d$ for takeoff equal or greater than the rate $\mu^d_{cong}$, which implies a service rate $\mu^a$ for landing equal or lower than the rate $\Psi(\mu^d_{cong})$.

A sustainable policy consists then on a policy $(\mu^a,\mu^d)$ within Gilbo's convex domain, such that $\mu^a\geq \mu^a_{cong}$ and $\mu^d\geq \mu^d_{cong}$.

Observe that Gilbo's sustainable policy domain will be empty if $\Psi(\mu^d_{cong})\leq \mu^a_{cong}$ or if $\Phi(\mu^a_{cong})\leq \mu^d_{cong}$.
A necessary and sufficient condition for the existence of a sustainable service policy is then:

	\begin{equation}\label{possivel}
		\mu^a_{cong}<\Psi(\mu^d_{cong}),\text{ equivalently } \mu^d_{cong}<\Phi(\mu^a_{cong})
	\end{equation}
Both conditions are equivalent because $\Phi,\Psi$ are monotone decreasing, inverse to each other.

We observe now that policies $(\mu^a,\mu^d)$ can be represented using service rates for the landing and takeoff operations but if the slot has a given demand rate $(\lambda^a,\lambda^d)$ then stable policies (those that lead to stable queues) can also be represented by the associated stable transit times $(z^a,z^d)$. We have a transformation from stable transit times in the positive quadrant $]0,+\infty[\times]0,+\infty[$ to service rates in the region  $]\lambda^a,+\infty[\times ]\lambda^d,+\infty[$ using the following transformations (inverse to each other):
	\begin{equation}\label{trocazmu}
		\begin{aligned}
		&(z^a,z^d)\mapsto (\mu(\lambda^a,z^a,q^a),\mu(\lambda^d,z^d,q^d))	\\	&(\mu^a,\mu^d)\mapsto (z(\lambda^a,\mu^a,q^a),z(\lambda^d,\mu^d,q^d))
		\end{aligned}
	\end{equation}
where $\mu(\lambda,z,q)$ and $z(\lambda,\mu,q)$ are given in (\ref{relatezmulambda}). Recall that for fixed $\lambda,q$ the functions $z(\lambda,\mu,q)$ and $\mu(\lambda,z,q)$ are monotone decreasing with limit values $z\to+\infty$ when $\mu\to \lambda^+$ and $z\to 0$ when $\mu\to +\infty$.

Any choice of point $(z^a,z^d)\in]0,+\infty[\times]0,+\infty[$ shall be called a delay policy. For a given configuration and depending on the demand rates, certain delay policies can be achieved using an appropriate policy within Gilbo's capacity domain associated to that configuration, and others may not.

The presence of a configuration, demand rates and predefined service levels for a given slot establishes a correspondence of the given configuration's Gilbo's capacity domain  and a corresponding domain of delay policies. Policies that lead to congestion or to saturation lie above a certain (unbounded) curve, the envelope of stable delay policies, which is the image of an arc from Gilbo's capacity envelope, as is represented in figure \ref{DelayPolicyDomain}.

\pgfkeys{/tikz/.cd,
	lamba/.store in=\lamba,
	lambd/.store in=\lambd,
	qa/.store in=\qa,
	qd/.store in=\qd,
	pd/.store in=\pd,
	pa/.store in=\pa
}
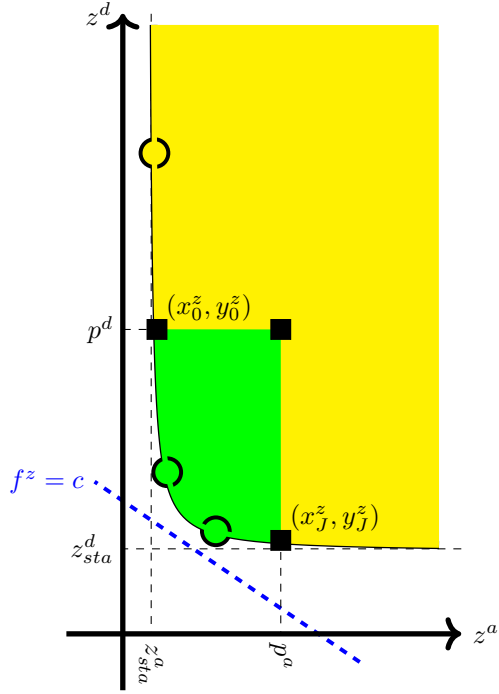
\begin{figure}
	\begin{tikzpicture}[lamba=2.,lambd=1.2, qa=2.1, qd=4.2, pa=1.4, pd=2.7, scale=1.5, declare function = {
			z(\la,\mu,\q) = 1/\mu + (\q*\la)/(2*\mu*(\mu-\la));
			la(\mu,\z,\q) = \mu*2*(\z*\mu-1)/(\q+2*(\z*\mu-1));
			mu(\z,\la,\q)=0.5*(1+\la*\z+(1+\la^2*\z^2+2*\la*\z*(\q-1))^0.5)/\z;
			valory(\valorx)=0.5*(6.5^2-\valorx^2)^0.5;
			valorx(\valory)=(6.5^2-4*\valory^2)^0.5;
		}
		]

		\draw[blue, line width=1.5pt, dashed] ({1.5*\pa},-0.25) -- (-0.25,{0.5*\pd}) node[left] {$f^z=c$};

		\draw[line width=2pt, ->] (-0.5,0) -- (3,0) node[right] {$z^a$};
		\draw[line width=2pt, ->] (0,-0.5) -- (0,5.5) node[left] {$z^d$};

		\draw[line width=3pt, variable=\t,domain={acos((mu(2*\pa,\lamba,\qa))/6.5):asin((mu(2*\pd,\lambd,\qd))/3.25)},samples=50]	plot ({z(\lamba,6.5*cos(\t),\qa)},{z(\lambd,3.25*sin(\t),\qd)});

		\draw[dashed]  ({z(\lamba,valorx(\lambd),\qa)},5.5)-- ({z(\lamba,valorx(\lambd),\qa)},0) node[right, rotate=-90]{$z^a_{sta}$};
		
		\draw[dashed] (0,{\pd}) node[left]{$p^d$}-- ({\pa},{\pd})-- ({\pa},0) node[right, rotate=-90]{$p^a$};;
		
		\draw[dashed] (3,{z(\lambd,valory(\lamba),\qd)})--(0,{z(\lambd,valory(\lamba),\qd)}) node[left] {$z^d_{sta}$};

		\fill [yellow, opacity=0.5, domain={acos((mu(2*\pa,\lamba,\qa))/6.5):asin((mu(2*\pd,\lambd,\qd))/3.25)}, variable=\t]
		({2*\pa}, {2*\pd})
		-- plot ({z(\lamba,6.5*cos(\t),\qa)},{z(\lambd,3.25*sin(\t),\qd)})
		-- ({2*\pa}, {2*\pd});

		\fill [green, opacity=0.5, domain={acos((mu(\pa,\lamba,\qa))/6.5):asin((mu(\pd,\lambd,\qd))/3.25)}, variable=\t]
		({\pa}, {\pd})
		-- plot ({z(\lamba,6.5*cos(\t),\qa)},{z(\lambd,3.25*sin(\t),\qd)})
		-- cycle;

		\foreach \t in {30,45,60} {
			\fill[color=black] ({z(\lamba,6.5*cos(\t),\qa)},{z(\lambd,3.25*sin(\t),\qd)}) circle [radius=4pt];
		}
		\foreach \t in {45,60} {
			\fill[color=green] ({z(\lamba,6.5*cos(\t),\qa)},{z(\lambd,3.25*sin(\t),\qd)}) 	circle [radius=3pt];}
		\fill[color=yellow] ({z(\lamba,6.5*cos(30),\qa)},{z(\lambd,3.25*sin(30),\qd)}) circle [radius=3pt];
		
		\draw[color=yellow,line width=2pt, variable=\t, domain={acos((mu(\pa,\lamba,\qa))/6.5)}:{acos((mu(2*\pa,\lamba,\qa))/6.5)},samples=50]
		plot ({z(\lamba,6.5*cos(\t),\qa)},{z(\lambd,3.25*sin(\t),\qd)});
		
		\draw[color=yellow,line width=2pt, variable=\t, domain={asin((mu(2*\pd,\lambd,\qd)/3.25)}:{asin((mu(\pd,\lambd,\qd))/3.25)},samples=50]
		plot ({z(\lamba,6.5*cos(\t),\qa)},{z(\lambd,3.25*sin(\t),\qd)});

		\draw[color=green,line width=2pt, variable=\t, domain={acos((mu(\pa,\lamba,\qa))/6.5)}:{asin((mu(\pd,\lambd,\qd))/3.25)},samples=50]
		plot ({z(\lamba,6.5*cos(\t),\qa)},{z(\lambd,3.25*sin(\t),\qd)});
		
		\draw [black] plot [only marks, mark size=2.5, mark=square*] coordinates {({\pa},{z(\lambd,valory(mu(\pa,\lamba,\qa)),\qd)}) } node[above right] {$(x^z_J,y^z_J)$};

		\draw [black] plot [only marks, mark size=2.5, mark=square*] coordinates {({z(\lamba,valorx(mu(\pd,\lambd,\qd)),\qa)},{\pd})} node[above right] {$(x^z_0,y^z_0)$};

		\draw [black] plot [only marks, mark size=2.5, mark=square*] coordinates {({\pa},{\pd})};


		
	\end{tikzpicture}
	\caption{Delay policies domain corresponding to situation of figure \ref{figuraGilboDecomposto}. Splitting of the domain between sustainable and congestion regions. Control points and linear cost function $f^z$ described in section \ref{sec5}.}\label{DelayPolicyDomain}
\end{figure}

When Gilbo's domain is polygonal with vertex at $(0,0)$ and additional vertices given as control points (\ref{ordercontrolpoints}) as represented in figure \ref{controlpoints}, it is determined by linear inequalities (\ref{segments}). At any slot with demand $(\lambda^a,\lambda^d)$ some of these control points will represent saturation conditions, some may represent congestion conditions and some fewer (if any) may represent sustainable conditions. To the effect of computing possible sustainable delay policies $(z^a,z^d)$ we should first determine sustainable service policies which are represented by (\ref{segments}) together with the additional conditions $\mu^a\geq \mu^a_{cong}$ and $\mu^d\geq \mu^d_{cong}$ using the congestion parameters given in (\ref{mucong}).

From a given configuration (hence a given Gilbo's convex domain), given predefined service levels $(p^a,p^d)$ and given demand rates $(\lambda^a,\lambda^d)$ we call sustainable policy any point $(\mu^a,\mu^d)$ of Gilbo's domain such that $\mu^a\geq \mu^a_{cong}$, $\mu^d\geq \mu^d_{cong}$. These points form a new convex domain whose envelope contains all the points of the arc $y=\Phi(x)$ with extremes at $(\mu^a_{cong},\Phi(\mu^a_{cong}))$ and $(\Psi(\mu^d_{cong}),\mu^d_{cong})$, and also the point $(\mu^a_{cong},\mu^d_{cong})$ and the vertical and horizontal line segments that join this point with the extreme points of the arc.

As Gilbo's domain is convex, assuming its envelope is a polygonal line with vertices at the origin and control points (\ref{ordercontrolpoints}), then Gilbo's sustainable policy domain will have as envelope a polygonal line determined by $(\mu^a_{cong},\mu^d_{cong})$ given by (\ref{mucong}) and alternative control points $((\bar x_j,\bar y_j))_{j=0\ldots J}$ given by:
	\begin{itemize}
		\item Control points $(\bar x_j,\bar y_j)=(x_j,y_j)$ when  $x_j\geq \mu^a_{cong}$ and $y_j\geq \mu^d_{cong}$.
		\item Coincident control points $(\bar x_j,\bar y_j)=(\mu^a_{cong},\Phi(\mu^a_{cong}))$ when $x_j<\mu^a_{cong}$, where $\Phi(\mu^a_{cong})$ is the maximal value of $\mu^d$ under the restrictions 	$$\mu^a_{cong}\cdot (y_j-y_{j-1})+\mu^d\cdot (x_{j-1}-x_j) \leq x_{j-1}\cdot y_j-x_j\cdot y_{j-1},\quad j=1\ldots J$$
		that is
		$$\Phi(\mu^a_{cong})=\min\limits_{j=1\ldots J}  
		\frac{x_{j-1}\cdot y_j-x_j\cdot y_{j-1}-\mu^a_{cong}\cdot (y_j-y_{j-1})}{x_{j-1}-x_j}$$
		\item Coincident control points $(\bar x_j,\bar y_j)=(\Psi(\mu^d_{cong}),\mu^d_{cong})$ when $y_j<\mu^d_{cong}$, where $\Psi(\mu^d_{cong})$ is the maximal value of $\mu^a$ under the restrictions 	$$\mu^a\cdot (y_j-y_{j-1})+\mu^d_{cong}\cdot (x_{j-1}-x_j) \leq x_{j-1}\cdot y_j-x_j\cdot y_{j-1}$$
		that is
		$$\Psi(\mu^d_{cong})=\min\limits_{j=1\ldots J}  
		\frac{x_{j-1}\cdot y_j-x_j\cdot y_{j-1}-\mu^d_{cong}\cdot (x_{j-1}-x_j)}{(y_j-y_{j-1})}$$
		\end{itemize}

This introduces new control points seen by a squared mark in figure \ref{figuraGilboDecomposto}. Gilbo's sustainable policy domain is then given by constraints:
	\begin{equation}\label{sustainablesegments}\begin{aligned}
		&\mu^a\geq \mu^a_{cong},\quad \mu^d\geq \mu^d_{cong}\\ 
		&\mu^a\cdot (\bar y_j-\bar y_{j-1})+\mu^d\cdot (\bar x_{j-1}-\bar x_j) &\leq \bar x_{j-1}\cdot \bar y_j-\bar x_j\cdot \bar y_{j-1},\quad j=1\ldots J
	\end{aligned}
	\end{equation}
In fact for the case of consecutive coincident control points, some of these restrictions are trivial $0\leq 0$ conditions.

Observe that if a point $(x_j,y_j)$ such that $x_j< \mu^a_{cong}$ and $y_j< \mu^d_{cong}$ existed on Gilbo's envelope, due to monotony of $\Phi$ we would have that no point of Gilbo's domain would belong to the region in $[\lambda^a,+\infty[\times [\lambda^d,+\infty[$, in which case all transit times associated to some service rate $(\mu^a,\mu^d)$ of Gilbo's domain would lead to dealys greater than the predefined levels $(p^a,p^d)$. This case should not appear as long as the scheduled flights $(\lambda^a,\lambda^d)$ are compatible with the predefined level of service (one should not schedule a slot with a demand rate that turns the desired level of service impossible to reach).

\section{Runway system performance optimization}\label{sec5}

Consider a slot with runway system configuration characterized by a certain Gilbo capacity envelope $y=\Phi(x)$ (equivalently $x=\Psi(y)$) and that any service policy $(\mu^a,\mu^d)$ in Gilbo's convex domain determines independent random variables $S^a, S^d$ for each service time of landing ($S^a$) and takeoff ($S^d$) operations at the slot, with expected values $1/\mu^a$ and $1/\mu^d$ and with quadratic ratios of momenta $q^a_S$, $q^d_S$ respectively (in particular both are independent of the service rates and greater that 1).

Consider that a schedule $(\lambda^a,\lambda^d)$ is know for the slot and that it determines independent random variables $C^a, C^d$ for inter-arrival times of consecutive clients of the landing ($C^a$) and takeoff ($C^d$) operations at the slot, with expected values $1/\lambda^a$ and $1/\lambda^d$ and with quadratic ratios of momenta $q^a_C$, $q^d_C$ respectively (in particular both are independent of the demand rates, are greater or equal than 1, and have value 2, if client arrival is a Poisson process).

Admit all general assumptions declared in section \ref{sec3}. In particular assume that $q^a=q^a_C+q^a_S-2\geq 1$ (which is the case in several situations, in particular if the clients arrive following a Poisson process). Assume the same property for takeoff operations: $q^d=q^d_C+a^d_S-2\geq 1$. Assume finally the estimation (\ref{Kingmanestimation}) for the expected queue lengths \cite{Kingman} holds in both queue systems (landing and takeoffs).

Transit times for airplanes demanding a landing service or a takeoff service is random with characteristics that evolve with time, and with expected value in the long term or in the stable situation called stable transit times $z^a,z^d$, given by formula (\ref{calculaz}), when $\lambda^a\mu^a$ and $\lambda^d<\mu^d$. 

Consider the airport operates with predefined delay tolerances $(p^a,p^d)$. We assume that if service policies have always been chosen in the sustainable region associated to this level of service, then the landing and takeoff queues at the beginning of the present slot are close to the stationary case corresponding to the parameters of the slot and hence assume that all clients arriving in this slot will have transit times of $z^a$ (in the case of landings) and of $z^d$ (in the case of takeoffs) computed using formula (\ref{calculaz}), with the appropriate parameters:
\begin{equation}\label{calculazazd}
	z^a=\frac1{\mu^a}\cdot \left(1+\frac{q^a\cdot \lambda^a}{2(\mu^a-\lambda^a)}\right),\qquad 
	z^d=\frac1{\mu^d}\cdot \left(1+\frac{q^d\cdot \lambda^d}{2(\mu^d-\lambda^d)}\right) 
\end{equation}
Recall that a necessary and sufficient condition for the existence of a sustainable service policy is any of the following:
$$\mu^a_{cong}<\Psi(\mu^d_{cong}),\qquad \mu^d_{cong}<\Phi(\mu^a_{cong})$$
with the congestion service ratios $\mu^a_{cong},\mu^d_{cong}$ given in formula (\ref{mucong}).

The expected cost due to delays in landing of flights scheduled at the present slot is usually a linear function on the aggregated transit times of all these flights, hence an expression $c^a\cdot \lambda^a\cdot z(\lambda^a,\mu^a,q^a)$, where $c^a$ is a coefficient that measures the cost associated to 1 slot of delay for 1 landing flight. In the same manner the expected cost due to delays in takeoff of flights scheduled at the present slot has an analogous form $c^d\cdot \lambda^d\cdot z(\lambda^d,\mu^d,q^d)$. The optimization problem to determine a policy of runway system with maximal performance will be the following:
	\begin{equation*}
	\begin{aligned}
		&\min f=c^a\cdot \lambda^a\cdot z(\lambda^a,\mu^a,q^a)+c^d\cdot \lambda^d\cdot z(\lambda^d,\mu^d,q^d) \\
		&\begin{aligned} \text{s.t. }  &\mu^d\leq \Phi(\mu^a) \\ &\mu^a\geq \mu^a_{cong},\, \mu^d\geq \mu^d_{cong} \end{aligned}
	\end{aligned}
\end{equation*}
which we call the runway system performance optimization problem for the current slot. It is characterized by the configuration (represented by $\Phi$, and the quadratic ratios of momenta $q^a,q^d$), and by the schedule (represented by demand ratios $(\lambda^a,\lambda^d)$) and where $\mu^a_{cong}$, $\mu^d_{cong}$, $z(\lambda,\mu,q)$ are given by (\ref{mucong}) and (\ref{relatezmulambda}). 

This is a mathematical program with decision variables $(\mu^a,\mu^d)$, a non-linear constraint given by the non-linear function $\Phi$ and a nonlinear objective function expressed in terms of the function $z(\lambda,\mu,q)$ from (\ref{calculaz}).

If we don't know an analytical specific expression for $\Phi(\mu)$ it is natural to use control points $((x_j,y_j))_{j=0\ldots J}$ as described in (\ref{ordercontrolpoints}), and to substitute the non-linear constraint described by $\Phi$ for the set of linear constraints given in (\ref{segments}), or the equivalent ones given in (\ref{sustainablesegments}), which only use control points that are in the stable region of the Gilbo envelope. 

However, this linearization of the constraints still leaves the problem with a non-linear objective function.

As the region of admissible points lies in the set of stable point of Gilbo's domain, we may apply transformations (\ref{trocazmu}), use the monotony properties of these functions, and express the optimization problem in terms of stable transit times $(z^a,z^d)$:
	\begin{equation*}
	\begin{aligned}
		&\min f^z=c^a\cdot \lambda^a\cdot z^a+c^d\cdot \lambda^d\cdot z^d \\
		&\begin{aligned} \text{s.t. }  &z^d\geq z(\lambda^d, \Phi(\mu(\lambda^a,z^a,q^a)),q^d) \\ &z^a\leq p^a,\, z^d\leq p^d \end{aligned}
	\end{aligned}
\end{equation*}
The curve $y=z(\lambda^d, \Phi(\mu(\lambda^a,x,q^a)),q^d)$ which may be interpreted as a sustainable transit time envelope, is not easy to represent but, as illustrated in figure \ref{DelayPolicyDomain}, it may be substituted by the polygonal line passing through transit time control points: $$(x^z_j,y^z_j)=(z(\lambda^a,\bar x_j,q^a), z(\lambda^d,\bar y_j,q^d)).$$ Adopting these points as belonging to the boundary of the sustainable transit times domain, the original problem is represented by a linear objective function $f^z(z^a,z^d)$, and constraints that are linearized as:
	\begin{equation*}\begin{aligned}
		&z^a\leq p^a,\quad z^d\leq p^d\\ 
		&z^a\cdot (y^z_j-y^z_{j-1})+z^d\cdot (x^z_{j-1}-x^z_j) &\leq x^z_{j-1}\cdot y^z_j-x^z_j\cdot y^z_{j-1},\quad j=1\ldots J
	\end{aligned}
\end{equation*}
Recall that the polygonal domain that describes sustainable delay policies uses as vertices $(x^z_j,y^z_j)$ the image by $(z(\lambda^a,x,q^a),z(\lambda^a,y,q^a))$ of sustainable service policies given as control points. We deduce that a basic optimal sustainable delay policy $(z^a,z^d)$ for this linear program is one of the vertices $(x^z_j,y^z_j)$ with associated cost $f$ not greater than for any of the remaining sustainable delay policies used as control points. Therefore one recovers sustainable service rates $\mu^a=\mu(\lambda^a,z^a,q^a)$ and $\mu^d=\mu(\lambda^d,z^d,q^d)$ which represent a sustainable service policy  $(\bar{x}_j,\bar{y}_j)$, where the cost $f$ is not greater than the cost associated to any of the remaining sustainable service rates given as control points.

To summarize: taking into account that Gilbo's convex domain is only known from empirical data, and that restricting condition $\mu^d\leq \Phi(\mu^a)$ for service policies depends on the concave function $\Phi$ which is not known, one may use the specific data $(\lambda^a,\lambda^d)$, $(p^a,p^d)$ associated to the present slot to translate Gilbo's domain into a new ``delay policies domain'' obtained by an appropriate transformation of the empirical data. Using a linear approximation to this new domain, and taking advantage of the linear expression of the performance function in these new variables, optimization of the runway system performance can be achieved by classical linear programming techniques.

\section{Sustainable demand rates domain}

When applying the previous performance optimization for the several time slots that an airport keeps on service each day, there arises a typical problem, namely condition (\ref{possivel}) may not hold at particular slots. In this case there does not exist a sustainable policy and the runway system performance optimization problem is meaningless. We might expect that some strategic measures would apply so that the flight schedule would conform to the capacities of the runway system. However, contracts to serve certain connections, or meteorological conditions that impose specific runway configurations might lead to a schedule that provokes congestion or saturated conditions  at certain slots. If such a situation arises, many of our assumptions in section \ref{sec5} do not hold anymore.

In this situation a natural technical solution is to transfer some of the scheduled flights from one slot to the following one. This imposes an additional slot of delay for all transferred flights, which would be added to the transit time that they will suffer until service is completed in the next slot. That is, at each time slot we must determine whether a flight slot transfer is needed.

Consider then a dynamic situation, of an airport with predefined delay tolerances $(p^a,p^d)$ for landing and takeoff operations (maximal transit times for airplanes demanding these operations, measured in slot size).  As stable transit times are given by (\ref{calculazazd}), the desired level of service $z^a\leq p^a$ and $z^d\leq p^d$ can be obtained only with service rates  $\mu^a\geq 1/p^a$ and $\mu^d\geq 1/p^d$.

Consider a continuous operation time, which is divided into $N$ consecutive time slots, that we identify by $i=1,\ldots, N$. A typical situation is a continuous operation time of 18 hours divided into $N=72$ consecutive time slots of 15min each. Consider that for each slot there is a schedule of $(\lambda^a_i,\lambda^d_i)$ landing and takeoff flights.

Consider that slot $i$ operates following a configuration determined by Gilbo's function, and quadratic ratio of momenta $(\Phi_i,q^a_i,q^d_i)$. Gilbo \cite{Airport_Capacity_Gilbo93} already dealt with this situation in the deterministic case, where the runway system could be managed with a demand rate equal to the service rate which, due to the deterministic nature of the model, caused no interference with the queue size. In Gilbo's paper, saturation was solved imposing a transfer of flights to the next slot, when needed to avoid saturation. As we have explained, if inter-arrival and service times are subject to whatever random factors, the service rate should stay a deal above the demand rate, in order to avoid congestion. Moreover, the cost of each decision is not linear on the service rates as assumed by Gilbo. It is not even quadratic on the service rate, as assumed by Jacquillat and Odoni, but would rather be linear in the stable transit times $z^a,z^d$, which have a non linear dependence on $\mu^a,\mu^d$. The optimization of this non linear cost component was studied in section \ref{sec5}. 

A natural question now is to determine which are the flight schedules $(\lambda^a,\lambda^d)$ that are compatible with sustainable service rates $(\mu^a,\mu^d)$ in Gilbo's domain. We know that each flight schedule $(\lambda^a,\lambda^d)$ determines (using formula (\ref{mucong})) specific minimal sustainable service rates $(\mu^a_{cong}(\lambda^a),\mu^d_{cong}(\lambda^d))\geq (1/p^a,1/p^d)$ and that services below these rates would lead to congestion.

Conversely, for any $(\mu^a,\mu^d)\geq (1/p^a,1/p^d)$ one may use formulas (\ref{relatezmulambda}) to determine the corresponding demand rates 	
	\begin{equation} \label{gilbosecond}
		\begin{aligned}
		\lambda^a_{cong}(\mu^a)&=\frac{2(p^a\mu^a-1)}{q^a+2(p^a\mu^a-1)}\cdot \mu^a\\
		\lambda^d_{cong}(\mu^d)&=\frac{2(p^d\mu^d-1)}{q^d+2(p^d\mu^d-1)}\cdot \mu^d		
		\end{aligned}
	\end{equation}
These demand rates are compatible with sustainable service rates $(\mu^a,\mu^d)$ and so are any demand rates $(\lambda^a,\lambda^d)\leq (\lambda^a_{cong},\lambda^d_{cong})$.  

We have then correspondences (\ref{mucong}) and (\ref{gilbosecond}) relating service rates and demand rates that are in correspondence with predefined stable transit times $z^a=p^a$, $z^d=p^d$ and therefore are limit cases of sustainable demand/service combinations for the predefined level of service. We also observe that functions $\mu^a_{cong}(\lambda)$, $\lambda^a_{cong}(\mu)$ are monotone increasing, inverse to each other, transforming the interval $[1/p^a,+\infty[$ into $[0,+\infty[$ and conversely. An analogous observation holds for  $\mu^d_{cong}(\lambda)$, $\lambda^d_{cong}(\mu)$.

We call sustainable demand rate domain, or briefly secondary Gilbo domain (See figure \ref{secondaryGilbo}) associated to a configuration represented by $(\Phi,q^a,q^d)$ and to specific levels of service $(p^a,p^d)$ the image using the transformation $(\lambda^a_{cong}(\mu^a),\lambda^d_{cong}(\mu^d))$ of Gilbo's domain $\mu^d\leq \Phi(\mu^a)$ on the region $\mu^a\geq 1/p^a$, $\mu^d\geq 1/p^d$. 

\pgfkeys{/tikz/.cd,
	lamba/.store in=\lamba,
	lambd/.store in=\lambd,
	qa/.store in=\qa,
	qd/.store in=\qd,
	pd/.store in=\pd,
	pa/.store in=\pa
}
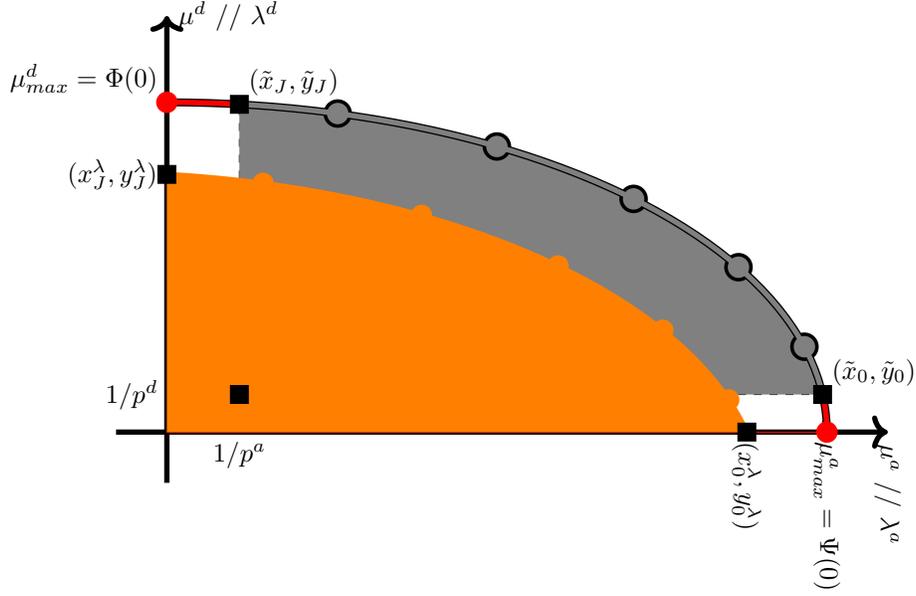
\begin{figure}
	\begin{tikzpicture}[lamba=0.,lambd=0, qa=2.1, qd=4.2, pa=1.4, pd=2.7, scale=1.35, declare function = {
			z(\la,\mu,\q) = 1/\mu + (\q*\la)/(2*\mu*(\mu-\la));
			la(\mu,\z,\q) = \mu*2*(\z*\mu-1)/(\q+2*(\z*\mu-1));
			mu(\z,\la,\q)=0.5*(1+\la*\z+(1+\la^2*\z^2+2*\la*\z*(\q-1))^0.5)/\z;
			valory(\valorx)=0.5*(6.5^2-\valorx^2)^0.5;
			valorx(\valory)=(6.5^2-4*\valory^2)^0.5;
		}
		]
		\draw[line width=2pt, ->] (-0.5,0) -- (7.1,0) node[right, rotate=-90] {$\mu^a$ // $\lambda^a$};
		\draw[line width=2pt, ->] (0,-0.5) -- (0,4.1) node[right] {$\mu^d$ // $\lambda^d$};
		\fill (0,3.25) circle (0.064cm) node[above left]{$\mu^d_{max}=\Phi(0)$};
		\fill (6.5,0) circle (0.064cm) node[right,rotate=-90]{$\mu^a_{max}=\Psi(0)$};
		\draw[dashed] (0,{mu(\pd,\lambd,\qd)}) node[left]{$1/p^d$}-- ({valorx(mu(\pd,\lambd,\qd))},{mu(\pd,\lambd,\qd)})-- ({valorx(mu(\pd,\lambd,\qd))},0);
		
		\draw[dashed] ({mu(\pa,\lamba,\qa)},0) node[below]{$1/p^a$}-- ({mu(\pa,\lamba,\qa)},{valory(mu(\pa,\lamba,\qa))})--(0,{valory(mu(\pa,\lamba,\qa))});
		
		
		\fill [red, opacity=0.5, domain=0:{asin(\lambd/3.25)}, variable=\t]
		({\lamba},0)-- plot ({6.5*cos(\t)},{3.25*sin(\t)})-- ({\lamba,\lambd})--cycle;	
		\fill [red, opacity=0.5, domain={acos(\lamba/6.5)}:90, variable=\t]
		(0,0)-- ({\lamba},0)-- plot ({6.5*cos(\t)},{3.25*sin(\t)})-- cycle;	
		
		
		\fill [gray, opacity=0.5, domain={acos((mu(\pa,\lamba,\qa))/6.5):asin((mu(\pd,\lambd,\qd))/3.25)}, variable=\t]
		({mu(\pa,\lamba,\qa)}, {mu(\pd,\lambd,\qd)})
		-- plot ({6.5*cos(\t)},{3.25*sin(\t)})
		-- cycle;	
		
		\fill [orange, opacity=0.5, domain={acos((mu(\pa,\lamba,\qa))/6.5):asin((mu(\pd,\lambd,\qd))/3.25)}, variable=\t]
		(0,0)
		-- plot ({la(6.5*cos(\t),\pa,\qa)},{la(3.25*sin(\t),\pd,\qd)})
		-- cycle;

		\draw[line width=3pt, variable=\t,domain={0:90},samples=50]	plot ({6.5*cos(\t)},{3.25*sin(\t))});
		
		\fill ({\lamba},{valory(\lamba)}) circle (0.0cm);
		\fill ({valorx(\lambd)},{\lambd}) circle (0.0cm);
		\foreach \t in {15,30,45,60,75} {
			\fill[color=black] ({6.5*cos(\t)},{3.25*sin(\t))}) circle [radius=4pt];
		}
		\foreach \t in {0,90} {
			\fill[color=red] ({6.5*cos(\t)},{3.25*sin(\t))}) circle [radius=3pt];
		}
		\foreach \t in {15,30,45,60,75} {
			\fill[color=gray] ({6.5*cos(\t)},{3.25*sin(\t))}) circle [radius=3pt];
		}
		\foreach \t in {15,30,45,60,75} {
		\fill[color=orange] ({la(6.5*cos(\t),\pa,\qa)},{la(3.25*sin(\t),\pd,\qd)}) circle [radius=3pt];
		}


		\draw[color=red,line width=2pt, variable=\t, domain={acos((mu(\pa,\lamba,\qa))/6.5)}:{acos(\lamba/6.5)},samples=50]
		plot ({6.5*cos(\t)},{3.25*sin(\t)});
		\draw[color=red,line width=2pt, variable=\t, domain={asin(\lambd/3.25)}:{asin((mu(\pd,\lambd,\qd))/3.25)},samples=50]
		plot ({6.5*cos(\t)},{3.25*sin(\t)});

		\draw[color=gray,line width=2pt, variable=\t, domain={acos((mu(\pa,\lamba,\qa))/6.5)}:{asin((mu(\pd,\lambd,\qd))/3.25)},samples=50]
		plot ({6.5*cos(\t)},{3.25*sin(\t)});

		\draw[color=orange,line width=2pt, variable=\t, domain={acos((mu(\pa,\lamba,\qa))/6.5)}:{asin((mu(\pd,\lambd,\qd))/3.25)},samples=50]
		plot ({la(6.5*cos(\t),\pa,\qa)},{la(3.25*sin(\t),\pd,\qd)});

		
		

		
		\draw [black] plot [only marks, mark size=2.5, mark=square*] coordinates {({mu(\pa,\lamba,\qa)},{valory(mu(\pa,\lamba,\qa))})} node[above right] {$(\tilde x_J,\tilde y_J)$};
		\draw [black] plot [only marks, mark size=2.5, mark=square*] coordinates { ({valorx(mu(\pd,\lambd,\qd))},{mu(\pd,\lambd,\qd)}) } node[above right] {$(\tilde x_0,\tilde y_0)$};
		
		\draw [black] plot [only marks, mark size=2.5, mark=square*] coordinates {({la(mu(\pa,\lamba,\qa),\pa,\qa)},{la(valory(mu(\pa,\lamba,\qa)),\pd,\qd)})} node[left] {$(x^\lambda_J,y^\lambda_J)$};
		\draw [black] plot [only marks, mark size=2.5, mark=square*] coordinates { ({la(valorx(mu(\pd,\lambd,\qd)),\pa,\qa)},{la(mu(\pd,\lambd,\qd),\pd,\qd)}) } node[right, rotate=-90] {$(x^\lambda_0,y^\lambda_0)$};

		\draw [black] plot [only marks, mark size=2.5, mark=square*] coordinates { ({mu(\pa,\lamba,\qa)},{mu(\pd,\lambd,\qd)})};
		
	\end{tikzpicture}
	\caption{For given delay tolerances, each Gilbo domain determines a secondary Gilbo domain, where all demand rates $(\lambda^a,\lambda^d)$ are compatible with a sustainable service rate associated to the configuration.}\label{secondaryGilbo}
\end{figure}

Taking into account the monotony of functions $\lambda^d_{cong}(\mu)$ and $\mu^a_{cong}(\lambda)$ described in (\ref{mucong}) and (\ref{gilbosecond}), an analytical characterization of Gilbo's secondary domain would be
	\begin{equation*}
	\begin{aligned}
		&\lambda^d\leq \left(\lambda^d_{cong}\circ \Phi\circ \mu^a_{cong}\right)(\lambda^a) \\ 
		&\lambda^a\geq 0,\, \lambda^d\geq 0
	\end{aligned}
	\end{equation*}

When Gilbo's domain is polygonal with vertex at $(0,0)$ and additional vertices given as control points (\ref{ordercontrolpoints}) as represented in figure \ref{controlpoints}, it is determined by linear inequalities (\ref{segments}). 

We pretend to transform points of Gilbo's domain that satisfy $x\geq 1/p^a$ and $y\geq 1/p^d$  using the mapping $(\lambda^a_{cong}(\mu^a),\lambda^d_{cong}(\mu^d))$, which is non-linear. This leads to a Gilbo secondary domain that is not polygonal. However, the assumption that the original Gilbo domain was polygonal was only assumed for convenience. What we really know is that a series of control points $(x_j,y_j)$ do belong to its envelope curve $y=\Phi(x)$. 

With similar arguments as applied in section \ref{sec4} consider new control points $(\tilde x_j,\tilde y_j)$ that belong to the region in $[1/p^a,+\infty[\times [1/p^d,+\infty[$  of Gilbo's envelope:
	\begin{itemize}
	\item Control points $(\tilde x_j,\tilde y_j)=(x_j,y_j)$ when  $x_j\geq 1/p^a$ and $y_j\geq 1/p^d$.
	\item Coincident control points $(\tilde x_j,\tilde y_j)=(1/p^a,\Phi(1/p^a))$ when $x_j<1/p^a$, where $\Phi(1/p^a)$ is the maximal value of $\mu^d$ under the restrictions 	$$(1/p^a)\cdot (y_j-y_{j-1})+\mu^d\cdot (x_{j-1}-x_j) \leq x_{j-1}\cdot y_j-x_j\cdot y_{j-1},\quad j=1\ldots J$$
	that is
	$$\Phi(1/p^a)=\min\limits_{j=1\ldots J}  
	\frac{x_{j-1}\cdot y_j-x_j\cdot y_{j-1}-(1/p^a)\cdot (y_j-y_{j-1})}{x_{j-1}-x_j}$$
	\item Coincident control points $(\tilde x_j,\tilde y_j)=(\Psi(1/p^d),1/p^d)$ when $y_j<1/p^d$, where $\Psi(1/p^d)$ is the maximal value of $\mu^a$ under the restrictions 	$$\mu^a\cdot (y_j-y_{j-1})+(1/p^d)\cdot (x_{j-1}-x_j) \leq x_{j-1}\cdot y_j-x_j\cdot y_{j-1}$$
	that is
	$$\Psi(1/p^d)=\min\limits_{j=1\ldots J}  
	\frac{x_{j-1}\cdot y_j-x_j\cdot y_{j-1}-(1/p^d)\cdot (x_{j-1}-x_j)}{(y_j-y_{j-1})}$$
\end{itemize} 
Observe that if a point $(x_j,y_j)$ such that $x_j<1/p^a$ and $y_j<1/p^d$ existed on Gilbo's envelope, due to monotony of $\Phi$ we would have that no point of Gilbo's domain would belong to the region in $[1/p^a,+\infty[\times [1/p^d,+\infty[$, in which case all transit times associated to some service rate $(\mu^a,\mu^d)$ of Gilbo's domain would lead to dealys greater than the predefined levels $(p^a,p^d)$. This case should not appear as long as the predefined levels $(p^a,p^d)$ are compatible with the properties of the runway system (one should not choose a lower bound for the service level of airplanes that is impossible to reach for the runway system).

The curve $y=(\lambda^d_{cong}\circ\Phi\circ\mu^a_{cong})(x)$ which bounds all demand rates compatible with a sustainable service, as illustrated in figure \ref{secondaryGilbo}, may be substituted by the polygonal line passing through demand rate control points: $$(x^\lambda_j,y^\lambda_j)=(\lambda^a_{cong}(\tilde x_j), \lambda^d_{cong}(\tilde y_j)).$$ Adopting these points as belonging to the boundary of the secondary Gilbo domain, our constraints are linearized as:
\begin{equation}\label{lambdasegments}\begin{aligned}
		&\lambda^a\geq 0,\quad \lambda^d\geq 0\\ 
		&\lambda^a\cdot (y^\lambda_j-y^\lambda_{j-1})+\lambda^d\cdot (x^\lambda_{j-1}-x^\lambda_j) &\leq x^\lambda_{j-1}\cdot y^\lambda_j-x^\lambda_j\cdot y^\lambda_{j-1},\quad j=1\ldots J
	\end{aligned}
\end{equation}
Recall that these control points depend on the runway configuration (meaning Gilbo's capacity envelope, and quadratic ratios of momenta for the services) and also on the specific service level (maximal admissible delays) that are active at a given slot.

\section{Optimization of flight slot transfers}

For a flight schedule $(\lambda^a_i,\lambda^d_i)_{i=1\ldots N}$ known for all the slots of a given day, we shall consider the problem to determine flight slot transfers. By this we mean to determine values $(s^a_i,s^d_i)$ (which we call a delay decision) representing a number of flights from slot $i$ to be transferred to the next slot, both for the landing and takeoff services. We impose $s^a_N=s^d_N=0$, that is, no transfer is admissible for the last slot of the day. We also impose $s^a_i,s^d_i\geq 0$. As our model is probabilistic in nature, there is no need to impose that these variables are integer, and a non integer number of airplanes to transfer can be seen as average number of transferred flights.

Consideration of $(s^a_i,s^d_i)_{i=1\ldots N}$ leads to a secondary flight schedule:
$$\lambda^{a2}_i=\lambda^a_i+s^a_{i-1}-s^a_i,\quad  \lambda^{d2}_i=\lambda^d_i+s^d_{i-1}-s^d_i$$
where we assume $s^a_0=s^d_0=0$.

\pgfkeys{/tikz/.cd,
	lamba/.store in=\lamba,
	lambd/.store in=\lambd,
	qa/.store in=\qa,
	qd/.store in=\qd,
	pd/.store in=\pd,
	pa/.store in=\pa
}
\begin{figure}
	\begin{tikzpicture}[lamba=2.,lambd=1.2, qa=2.1, qd=4.2, pa=1.4, pd=2.7, scale=1.35, declare function = {
			z(\la,\mu,\q) = 1/\mu + (\q*\la)/(2*\mu*(\mu-\la));
			la(\mu,\z,\q) = \mu*2*(\z*\mu-1)/(\q+2*(\z*\mu-1));
			mu(\z,\la,\q)=0.5*(1+\la*\z+(1+\la^2*\z^2+2*\la*\z*(\q-1))^0.5)/\z;
			valory(\valorx)=0.5*(6.5^2-\valorx^2)^0.5;
			valorx(\valory)=(6.5^2-4*\valory^2)^0.5;
		}
		]
		\draw[line width=2pt, ->] (-0.5,0) -- (7.1,0) node[right] {$\mu^a$};
		\draw[line width=2pt, ->] (0,-0.5) -- (0,4.1) node[left] {$\mu^d$};
		\draw[dashed] (0,{\lambd}) node[left]{$\lambda^{d}_i+s^d_{i-1}-s^d_i$}-- ({valorx(\lambd)},{\lambd})-- ({valorx(\lambd)},0);
		\draw[dashed] (0,{mu(\pd,\lambd,\qd)}) -- ({valorx(mu(\pd,\lambd,\qd))},{mu(\pd,\lambd,\qd)})-- ({valorx(mu(\pd,\lambd,\qd))},0);
		
		\draw[dashed] ({\lamba},0) node[below]{$\lambda^{a}_i+s^a_{i-1}-s^a_i$}-- ({\lamba},{valory(\lamba)})--(0,{valory(\lamba)});
		\draw[dashed] ({mu(\pa,\lamba,\qa)},0)-- ({mu(\pa,\lamba,\qa)},{valory(mu(\pa,\lamba,\qa))})--(0,{valory(mu(\pa,\lamba,\qa))});
		
		\fill [green, opacity=0.5, domain={acos((mu(\pa,\lamba,\qa))/6.5):asin((mu(\pd,\lambd,\qd))/3.25)}, variable=\t]
		({mu(\pa,\lamba,\qa)}, {mu(\pd,\lambd,\qd)})
		-- plot ({6.5*cos(\t)},{3.25*sin(\t)})
		-- cycle;	
		
		\draw[line width=3pt, variable=\t,domain={0:90},samples=50]	plot ({6.5*cos(\t)},{3.25*sin(\t))});
		
		\foreach \t in {0,15,30,45,60,75,90} {
			\fill[color=black] ({6.5*cos(\t)},{3.25*sin(\t))}) circle [radius=4pt];
		}
		\foreach \t in {45,60} {
			\fill[color=green] ({6.5*cos(\t)},{3.25*sin(\t))}) circle [radius=3pt];
		}
		\draw[color=green,line width=2pt, variable=\t, domain={acos((mu(\pa,\lamba,\qa))/6.5)}:{asin((mu(\pd,\lambd,\qd))/3.25)},samples=50] plot ({6.5*cos(\t)},{3.25*sin(\t)});
		
		\fill[color=orange]{(\lamba,\lambd)} circle [radius=4pt];
		\fill[color=blue](5.5,1.4) circle [radius=4pt]  node[left]{$(\lambda^a_i+s^a_{i-1},\lambda^d_i+s^d_{i-1})$};
		\draw[line width=2pt, dashed, -{Stealth[length=3mm]}] (5.5,1.4) to[bend left] (2,1.2);
		
	\end{tikzpicture}
	\caption{For a given (scheduled plus previously delayed) demand rate at a slot, it may even lie on the non saturating domain, but if it is not compatible with a sustainable service, we explore the possibility to diminish the demand rate by transferring $(s^a_i,s^d_i)$ flights to the next slot.}\label{flightslottransfer}
\end{figure}
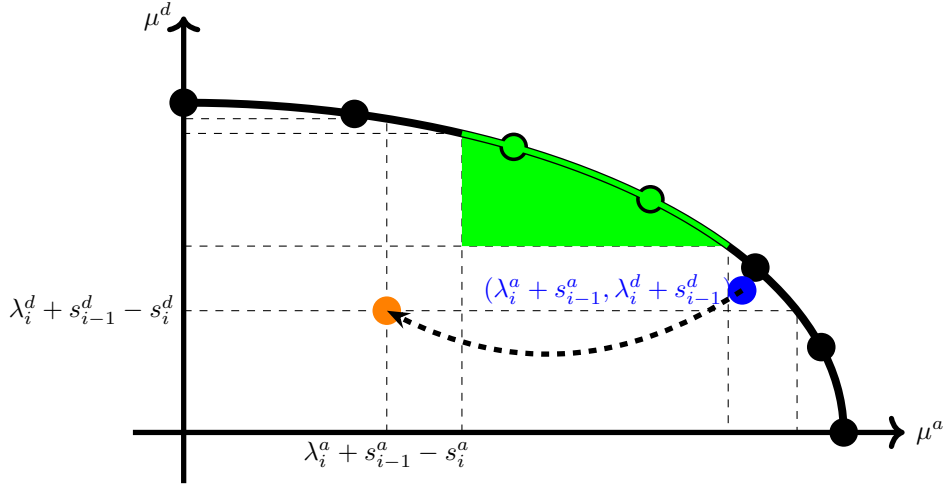

In order to guarantee an admissible delay decision, we have the additional constraints $\lambda^{a2}_i\geq 0$, $\lambda^{d2}_i\geq 0$.

Furthermore, we will impose that the secondary flight schedule $(\lambda^{a2}_i,\lambda^{d2}_i)$ lies, for each slot, in the secondary Gilbo domain (hence policies can be applied that give a sustainable service). 

We assume that one such delay choice implies a cost $c^a$ for each transfer of a landing airplane to the next slot and a cost $c^d$ for each transfer of a takeoff airplane to the next slot. The total cost of such a decision would be:
$$f=c^a\cdot \sum_i s^a_i+ c^d\cdot \sum_i s^d_i$$
We shall consider the following optimization problem that leads to a sustainable schedule with minimal cost of the delay decisions:

	\begin{equation*}
	\begin{aligned}
		&\min f=c^a\cdot \sum_i s^a_i+ c^d\cdot \sum_i s^d_i \\
		&\begin{aligned} \text{s.t. }  & s^a_i,s^d_i\geq 0\quad (i=0\ldots N) \\ &s^a_0=s^d_0=s^a_N=s^d_N=0 \\
			\\
			&(\ref{lambdasegments}) \text{ holds for }
			\left\{\begin{aligned} &\lambda^{a2}_i=\lambda^a_i+s^a_{i-1}-s^a_i\\ &\lambda^{d2}_i=\lambda^d_i+s^d_{i-1}-s^d_i\\
			&\text{control points }(x^\lambda_j,y^\lambda_j)\text{ associated to }\\ &\text{configuration }(\Phi_i,q^a_i,q^d_i) \end{aligned}\right\} \end{aligned}
	\end{aligned}
	\end{equation*}
Observe that to solve this mathematical program it is interesting to previously compute all demand rate control points $(x^\lambda_j,y^\lambda_j)$ on the secondary Gilbo envelopes associated to each configuration of the runway system and with the specific service level $(p^a,p^d)$ adopted by the airport. The original operational throughput envelope is not necessary if the service level is maintained constant in all slots.

After solving this mathematical problem we get a  secondary schedule $(\lambda^{a2},\lambda^{d2})$ that is admissible, in the sense that the demand rates represented by this schedule allow for the airport to operate in sustainable conditions, when the appropriate policies are adopted. This secondary schedule has minimal cost among them. 

For the slots $i$ where $(\lambda^{a2}_i,\lambda^{d2}_i)$ is a vertex $(x^{\lambda}_j,y^{\lambda}_j)$ of Gilbo's secondary domain, in fact the slot will be in a sustainable situation with extreme transit times $z^a=p^a$, $z^d=p^d$ achieved by a single possible policy $(\mu^a,\mu^d)=(\tilde x_j,\tilde y_j)$. For the slots $i$ where $(\lambda^{a2}_i,\lambda^{d2}_i)$ is not such a vertex, a runway system performance optimization can be applied (section \ref{sec5}) to further minimize the costs associated to (not extremal) stable transit times.

\section{Conclusions and further work} 

The most essential information used to characterize transit times for landing and takeoff airplanes served by a runway system is described in terms of an operational throughput envelope together with quadratic ratios of momenta for the inter-arrival and service times. The  operational throughput envelope is empirical by its nature and is then characterized by a finite family of control points, rather than any explicit analytical expression.

Such information allows to estimate the expected delays due to queuing in the stable case, and to determine how these delays relate to each other depending on the specific service priority equilibrium choice between landing or takeoff.  It is assumed that this relation folows Kingman's \cite{Kingman} estimation of queue size, but any other functional relation between transit times $z$, demand rate $\lambda$ and service rate $\mu$, whith the obvious monotonicity properties, would also apply for this kind of study. In particular, such relations as described by \cite{Kim}

Fixing a given delay tolerance for the runway system operations, the operational throughput domain determines a secondary domain that contains all operation demand rates that can be assumed by the runway system while operating in sustainable conditions. This domain can be used to determine a linear problem that leads from any flight schedule to a secondary flight schedule obtained by flight slot transfers and that avoids the congestion of the airport, while minimizing the costs associated to such transfers.

Moreover, given a flight schedule that avoids the congestion of the airport, a nonlinear program arises when trying to apply policies (equilibrium of landing versus takeoff priorities) that minimizes the aggregate costs associated to transit times of all airplanes. This nonlinear program is solved using a description of the policy in terms of intended delays instead of using intended service rates. With these new parameters, the objective function is linear and control points identify a linearization of the region corresponding to sustainable operation of the runway system.

Our approach doesn't assume any specific probabilistic model for arrival or service times, but assumes the validity of Kingman's extension of Pollakzec-Khintchine formula when the runway operates in sustainable conditions. It represents an alternative to deterministic approaches or other approaches in the literature that considered a cost expressed as weighed mean of the squared queue lengths along the day, with specific models (Exponential, Erlang) for all involved random variables.

The present approach determines a functional relation between arrival rate, service rate, and transit time. It suggests the likely utility of registering historical data containing these three components of information, for each slot of operation with a specific runway system configuration and to determine up to what point a small sustainable queue at the beginning of any slot leads to a significant discrepancy of the observed values and those predicted by Kingman's formula.

This approach also identifies the most basic information needed to model a runway system in such a way that its performance (in terms of transit times) can be estimated. Its simplicity allows its use to determine characteristics of the system for a given airport, for example the determination of marginal costs associated to delay tolerances or to parameters of the operational throughput envelope. Such studies might be interesting for strategical decisions like choosing between a new runway system project or an enhancement of airport services that might widen its delay tolerance.

The model might also be the skeleton for the design of more advanced tools to deal with the complexity of runway traffic control, which might include other variables. For example when considering service times in terms of airplane size, or when the existence of connecting flights is relevant. This basic structure could be extended for those cases, simplifying new decision processes related to landing/takeoff priorities. This would imply  a distinction of more than 2 types of services and a collection of corresponding historical data. The whole theory might be also applied in other situations where a single server is used for competing queues of clients with different needs, demand rates that vary with time, and where a manager can adjust the proportion of services to be allocated to these client classes, taking into account predefined service levels and client arrival schedules.

\end{document}